\documentclass[final]{amsart}

\usepackage{graphicx}
\usepackage{amssymb}
\usepackage{mathtools}
\usepackage{relsize,exscale}

\newtheorem{theorem}{Theorem}
\newtheorem{lemma}{Lemma}[section]
\newtheorem{prop}[lemma]{Proposition}
\newtheorem{coroll}[lemma]{Corollary}
\theoremstyle{definition}
\newtheorem{definition}[lemma]{Definition}
\theoremstyle{remark}
\newtheorem{remark}[lemma]{Remark}

\newcommand{\bl}{\begin{lemma}}
\newcommand{\el}{\end{lemma}}
\newcommand{\bp}{\begin{prop}}
\newcommand{\ep}{\end{prop}}
\newcommand{\bcor}{\begin{coroll}}
\newcommand{\ecor}{\end{coroll}}
\newcommand{\bpf}{\begin{proof}}
\newcommand{\epf}{\end{proof}}
\newcommand{\brem}{\begin{remark}}
\newcommand{\erem}{\end{remark}}
\newcommand{\eps}{\varepsilon}

\begin{document}

\title[The asymptotic number of $BCK(p)$-terms]{The asymptotic number of $BCK(p)$-terms}

\author{Olivier Bodini}
\address{LAGA, CNRS UMR 7539, Institut Galil\'ee, Universit\'e Sorbonne Paris Nord (Paris 13),
99 avenue Jean-Baptiste Cl\'ement, 93430 Villetaneuse, France}
\email{olivier.bodini@univ-paris13.fr}
\thanks{O.~Bodini was supported by the French Agence Nationale de la Recherche (ANR) under the
PANDAG project, reference ANR-23-CE48-0014.}

\begin{abstract}
A $BCK(p)$-term is a closed lambda-term in which each abstraction binds at most $p$ variables;
for $p=1$ these terms are related to the axiom system BCK of combinatory logic. We determine the
number of $BCK(p)$-terms of size $n$, up to a constant multiplicative factor, for every $p\ge2$.
For $p=2$ the result is fully explicit:
\[
f_n\sim A\,n^{2n/5}2^{n/5}e^{-2n/5}n^{-3/5}e^{P(n^{1/5})},
\]
where $P$ is an explicit polynomial of degree $4$.
In general the asymptotic formula involves $n^{pn/(2p+1)}$ and an exponential of a polynomial in
$n^{1/(2p+1)}$. The proof combines a linearization of the defining functional equation, the
analysis of the Borel transform of the linearized generating function, which satisfies a linear
differential equation with an irregular singularity at infinity, the saddle point method for
Hayman-admissible functions, and positivity arguments comparing the original and the linearized
sequences. For $p=1$, where the asymptotic number of terms was determined earlier by a
different method, this comparison fails: the quadratic term then changes the polynomial factor,
by a factor of order $n^{1/3}$.
\end{abstract}

\maketitle

\section{Introduction}\label{sec:intro}

A lambda-term is a formal expression built of variables and a quantifier $\lambda$ which in
general occurs more than once and acts on one of the free variables. It can be described by the
context-free grammar $T::= a\;\mid\;(T*T)\;\mid\;\lambda a.T$ where $a$ is a variable. The
concatenation $T_1*T_2$ of terms is called \emph{application} and putting the string $\lambda a$
in front of a term is called \emph{abstraction}. Each abstraction $\lambda a$ \emph{binds} the
variable $a$ in the subsequent term and each variable can only be bound by at most one
abstraction. A term where all the variables are bound is called a \emph{closed} lambda-term.

Lambda-terms are the words of the formal language known as lambda-calculus (see \cite{Ba84}),
which was introduced by Church in the 1930s, with fundamental contributions by Kleene (see
\cite{Ch36, Kl35a, Kl35b}), in order to describe computable functions; it is now used to analyze
programming languages or investigate decision problems. Moreover, it is the basis of the
programming language LISP. Furthermore, by the Curry-Howard isomorphism (see \cite{SU06}),
simply typed lambda-terms can also be seen as proofs of intuitionistic tautologies.

Our original motivation was to gain information on the typical shape of large closed
lambda-terms, for which the number of closed lambda-terms of a given size has to be known. This
turns out to be a very difficult problem (see \cite{BGG11} for a brief discussion). We therefore
study particular subsets $\mathcal B_p$ of the set $\mathcal L$ of closed lambda-terms, which we
call $BCK(p)$-terms, $p=1,2,\dots$: $\mathcal B_p$ is the set of closed lambda-terms in which each
abstraction binds at most $p$ variables. The name comes from the connection of $\mathcal B_1$ with
the axiom system BCK in combinatory logic (see \cite{KT78} for BCK logic and \cite{BGGJ13} for an
explanation of the connection). We have
$\mathcal B_1\subset\mathcal B_2\subset\dots\subset\mathcal L$ and $\mathcal B_\infty=\mathcal L$.

The enumeration of $BCK(p)$-terms is of interest in its own right, because several difficulties
have to be overcome. First, we are facing unlabelled objects, so the formal specification gives
rise to a functional equation for an \emph{ordinary} generating function; but the number of terms
grows so fast that its radius of convergence is zero. Second, the recurrence relation for the
counting sequence is not holonomic, since it involves a quadratic term. Omitting the quadratic
term, we get a holonomic sequence; but the differential equation satisfied by the Borel transform
of its generating function has an irregular singularity at infinity, of a type in which all the
eigenvalues of the leading matrix coincide.

\subsection*{Main results}
Let $f^{(p)}_n$ denote the number of $BCK(p)$-terms of size $n$ (the size is the total number of
variables, abstractions and applications, see Section~\ref{sec:prelim}), and write $f_n=f^{(2)}_n$.

\begin{theorem}\label{thm:main2}
The number of $BCK(2)$-terms of size $n$ satisfies, as $n\to\infty$,
\begin{multline}\label{asymBCK}
f_n\sim A\, n^{2n/5} 2^{n/5} e^{-2n/5}\, n^{-\frac{3}{5}}\\
\times\exp\left(2^{-8/5}n^{4/5}+\frac{7\cdot 2^{4/5}}{15}n^{3/5}-\frac{17\cdot 2^{1/5}}{75}n^{2/5}
-\frac{41\cdot 2^{3/5}}{500}n^{1/5}\right),
\end{multline}
where $A$ is a positive constant.
\end{theorem}

The constant $A$ is not known explicitly; the convergence is extremely slow, and Manuel Kauers
\cite{Kau13} obtained $A\approx0.9223889$ from the first $500000$ terms.

For general $p$ we set
\begin{gather*}
\mu_p=\frac p{2p+1},\qquad \nu_p=\frac{2p+1}{p+1},\qquad
\varrho_p=\Bigl(\frac{2^p}{p!}\Bigr)^{1/(2p+1)},\\
\tau_p=\frac{p+1}{2p+1}\Bigl(\frac{2^p}{p!}\Bigr)^{1/(p+1)}=(1-\mu_p)\varrho_p^{\nu_p}.
\end{gather*}
For $p=2$ these are $\mu_2=\frac25$, $\nu_2=\frac53$, $\varrho_2=2^{1/5}$ and
$\tau_2=\frac{3\cdot2^{1/3}}5$.

\begin{theorem}\label{thm:mainp}
\begin{enumerate}
\item[(i)] For every $p\ge1$, $\log f^{(p)}_n=\mu_p\log n!+n\log\varrho_p+o(n)$.
\item[(ii)] For every $p\ge2$ there are constants $A_p>0$, $\sigma_p\in\mathbb R$ and
$e_1,\dots,e_{2p}\in\mathbb R$ such that
\[
f^{(p)}_n\sim A_p\,n^{\mu_pn}\bigl(\varrho_pe^{-\mu_p}\bigr)^n
\exp\Bigl(\sum_{m=1}^{2p}e_mn^{m/(2p+1)}\Bigr)\,n^{\sigma_p/\nu_p}.
\]
\end{enumerate}
\end{theorem}

The constants $e_m$ and $\sigma_p$ can be computed by a triangular recursion; for $p=2$,
Theorem~\ref{thm:mainp}(ii) is Theorem~\ref{thm:main2}.

The case $p=1$ ($BCK$-terms, also called affine lambda-terms) was solved by Bodini, Gardy and
Jacquot \cite{BGJ10a}, by means of a bijection with combinatorial maps. Their result can be
written as
\[
f^{(1)}_n\sim C_1\,n^{-1/6}\Bigl(\frac{2n}{e}\Bigr)^{n/3}
\exp\Bigl(\tfrac12(2n)^{2/3}-\tfrac16(2n)^{1/3}\Bigr).
\]
This has the form of Theorem~\ref{thm:mainp}(ii), except that the power $n^{\sigma_1/\nu_1}=n^{-1/2}$
is replaced by $n^{-1/6}$. Combining this result with our analysis of the linearized sequence, we
show that $f^{(1)}_n/g^{(1)}_n\sim c\,n^{1/3}$ (Remark~\ref{rem:p1}). For $p=1$ the quadratic term
therefore changes the polynomial factor, which explains why our comparison argument is
restricted to $p\ge2$.

\subsection*{Method}
Removing the quadratic term from the functional equation for $F_p(z)=\sum_nf^{(p)}_nz^n$ yields a
linear equation, whose solution $G_p(z)=\sum_ng^{(p)}_nz^n$ is D-finite. Its Borel transform
$H_p(z)=\sum_ng^{(p)}_nz^n/n!$ is entire and satisfies a linear differential equation with an
irregular singular point at infinity. Using the Newton polygon of this equation, the
Hukuhara--Turrittin theorem and the positivity of the coefficients, we prove
(Theorem~\ref{thm:Hp}) that, in a sector around the positive real axis,
\[
H_p(z)=C\,z^{\sigma_p}\exp\Bigl(\tau_p\,z^{(2p+1)/(p+1)}+\sum_{m=1}^{2p}a_mz^{m/(p+1)}\Bigr)
\bigl(1+O(|z|^{-1/(p+1)})\bigr).
\]
We then show that $H_p$ is Hayman-admissible and apply the saddle point method. Finally, we prove
that $f^{(p)}_n/g^{(p)}_n$ tends to a constant for $p\ge2$. This step rests on a decomposition of
$f^{(p)}$ along the Green's function of the linear recurrence and on an a priori bound
$f^{(p)}_n\le MK^n(n!)^{\mu_p}$, obtained by a descent argument in $K$.

Exponentials of fractional powers are not common in combinatorial enumeration. The only example
known to us with similar phenomena is the asymptotic enumeration of certain Young tableaux
\cite{CMS05}; there, however, the generating function involves integer powers only.

\subsection*{Related work}
The interest in counting lambda-terms or terms of other logical systems has recently arisen in
the literature. E.g. in \cite{DZ09} proofs in propositional systems are counted. David et
al.~\cite{DGKRTZ13} study the number of terms using a different definition of size from that we
use in this paper and analyze terms w.r.t. properties like normalizability. Grygiel and Lescanne
\cite{GryLes13} aim at generation of terms and derive several recursive formulas related to the
number of lambda-terms. Bounds for the number of terms as well as asymptotics for the number of
particular terms have been derived in \cite{BGGJ13}. Further work where terms with certain
restrictions are counted are \cite{BGG11} or \cite{BGJ10a}, where the asymptotic number of terms
related to the axiom systems BCI and BCK (i.e. our set $\mathcal B_1$) has been determined. A
comparison of these two systems was also carried out in \cite{GIZ13}. A preliminary version of the
present work, devoted to the case $p=2$, appeared in \cite{BG14}.

\subsection*{Plan of the paper}
Section~\ref{sec:prelim} defines $BCK(p)$-terms, recalls the functional equation of
\cite{BGGJ13}, derives the linearized equation for every $p$, and collects the analytic tools. Sections~\ref{sec:lin2}
to~\ref{sec:f2} treat the case $p=2$ in detail.
\begin{itemize}
\item Section~\ref{sec:lin2} establishes the recurrence and the growth of the linearized sequence.
\item Section~\ref{sec:borel2} analyzes its Borel transform: differential equation, asymptotics in
a sector, Hayman admissibility and saddle point.
\item Section~\ref{sec:f2} compares the original and the linearized sequences and proves
Theorem~\ref{thm:main2}.
\end{itemize}
Section~\ref{sec:general} extends the method to every $p$ and proves Theorem~\ref{thm:mainp}.
Section~\ref{sec:concl} discusses open questions. The specifications and their translation into
equations for generating functions are based on the symbolic method \cite[Ch.~I--III]{FlSe09}.

\section{$BCK(p)$-terms and the linearized equation}\label{sec:prelim}

First of all, we have to specify precisely the objects we want to count: We count closed lambda
terms up to $\alpha$-conversion (see \cite[Ch.~2]{Ba84}), that is, only the relation between
variables and quantifiers are significant, not the variable names themselves. In this context, we
prefer another viewpoint to look at lambda-terms:

\begin{definition} \label{deflambda}
A \emph{$BCK(p)$-term} is a directed acyclic graph constructed as follows: Start with a Motzkin
tree, i.e. a rooted plane tree where each node has 0, 1 or 2 children. This tree is the
\emph{underlying tree} of a lambda-term. The nodes with 0 children are called
\emph{leaves},\footnote{Note the slight abuse of the terminology: We add further edges such that
the resulting graph does not possess any end nodes. Still, when speaking of lambda-terms we always
keep the underlying tree in mind, so that such a wording does make sense.} those having 1 child
\emph{unary nodes} and those with 2 children \emph{binary nodes} of the lambda-term. Assume that
all the edges are directed away from the root. Then add further directed edges, called
\emph{pointers}, each having a unary node as start vertex and as end vertex a leaf which belongs
to the subtree rooted at the start vertex, in such a way that every leaf is the end vertex of
exactly one pointer. Moreover, every unary node must be the start vertex of at most $p$ pointers.

The \emph{size} of a $BCK(p)$-term is defined as the total number of its vertices.
\end{definition}

\begin{remark}
The correspondence between terms and graphs is obvious: Binary nodes correspond to applications,
unary nodes to abstractions (which bind at most $p$ leaves) and leaves to variables. Pointers
represent the binding relation.
\end{remark}

\begin{remark}
The size of a lambda-term as defined above satisfies the following recursive relation:
\[
|x|=1, \quad |\lambda x.T|=1+|T|, \quad |(S*T)|=1+|S|+|T|.
\]
Here $x$ denotes a variable and $S,T$ denote subterms. That means that each of the variables,
abstractions and applications of a term contribute 1 to its size.
\end{remark}

Now, we will recall some known results. We are interested in the numbers $f_n^{(p)}$ of
$BCK(p)$-terms of size $n$ and in order to derive an asymptotic expression for $f_n^{(p)}$ we
will work with the generating function
\[
F_p(z)=\sum_{n\ge 0} f_n^{(p)}z^n.
\]
We write $f_n=f_n^{(2)}$ and $F(z)=F_2(z)$.

The following functional equation was established in \cite{BGGJ13}.

\bp[\cite{BGGJ13}]\label{BCKdirect}
Let $M(z,u)$ denote the generating function of Motzkin trees where $z$ marks the size (i.e. the
total number of nodes) and $u$ marks the number of leaves. This function is given by the unique
power series solution of $M(z,u)=uz+zM(z,u)+zM(z,u)^2$, that is
\begin{equation} \label{motz}
M(z,u)=\frac{1-z-\sqrt{(1-z)^2-4uz^2}}{2z}.
\end{equation}
Then $F_p(z)$ is given as the solution of
\begin{equation}\label{bck}
F_p(z)=z[u^p]\frac{M(z,u)}{1-u} + z F_p(z)^2 + z [u^p]\frac1{1-u}F_p\left(\frac{z}{1-2zM(z,u)}\right).
\end{equation}
\ep

\subsection{The linearized equation}
Because of the quadratic term $zF_p^2$, the coefficients of $F_p$ do not satisfy a linear
recurrence with polynomial coefficients. We therefore consider the solution $G_p(z)=\sum_n
g^{(p)}_nz^n$ of \eqref{bck} with the quadratic term removed:
\begin{equation}\label{linearizedeqp}
G_p(z)=z[u^p]\frac{M(z,u)}{1-u} + z [u^p]\frac1{1-u}G_p\left(\frac{z}{1-2zM(z,u)}\right).
\end{equation}
Since $zF_p^2$ has nonnegative coefficients, $f^{(p)}_n\ge g^{(p)}_n$ for all $n$. We will show
that, for $p\ge2$, the quotient $f^{(p)}_n/g^{(p)}_n$ tends to a constant. The next lemma turns
\eqref{bck} and \eqref{linearizedeqp} into differential equations.

\bl\label{lem:linp}
Let $w=z/(1-2zM(z,u))$, let $t=z/(1-z)$, and set
\[
\Gamma_k(z)=\frac{z}{k!}\sum_{l=0}^p[u^l](w-z)^k\quad(0\le k\le p),\qquad
\mathcal A_p(z)=z\sum_{l=1}^p[u^l]M(z,u).
\]
Then $G_p=\mathcal A_p+\sum_{k=0}^p\Gamma_kG_p^{(k)}$ and
$F_p=\mathcal A_p+\sum_{k=0}^p\Gamma_kF_p^{(k)}+zF_p^2$. Moreover:
\begin{enumerate}
\item[(i)] $\Gamma_0=z$; for $k\ge1$, $\Gamma_k$ has nonnegative coefficients and valuation
$3k+1$; and $\Gamma_p(z)=\frac{2^p}{p!}\,\frac{z^{3p+1}}{(1-z)^p}$.
\item[(ii)] $(1-z)^{2p-1}\Gamma_k(z)$ and $(1-z)^{2p-1}\mathcal A_p(z)$ are polynomials. As
$z\to\infty$, $\Gamma_k(z)=O(z^{k+p+1})$ and $\mathcal A_p(z)=O(z)$.
\end{enumerate}
\el

\bpf
We use the Taylor expansion $G(w)=\sum_kG^{(k)}(z)(w-z)^k/k!$ and
$[u^p]\frac1{1-u}X=\sum_{l\le p}[u^l]X$. The function $M$ satisfies $M=ut+tM^2$,
hence
\[
M=\sum_{l\ge1}C_{l-1}t^{2l-1}u^l,
\]
where $C_{l-1}$ are the Catalan numbers. Since $w-z=\sum_{i\ge1}z(2zM)^i$, expanding gives
\[
[u^l](w-z)^k=z^k\sum_{I=k}^{l}\gamma_{k,l,I}\,(2z)^I\,t^{2l-I}
\]
with nonnegative integers $\gamma_{k,l,I}$, and $[u^l](w-z)^k=0$ for $l<k$. All the claims
follow:
\begin{itemize}
\item the smallest power of $z$ is $z^{k+2l}\ge z^{3k}$;
\item the pole at $z=1$ has order $2l-I\le2p-1$;
\item $t\to-1$ as $z\to\infty$, which gives the growth bounds;
\item for $k=p$ only $l=I=p$ occurs, with $\gamma_{p,p,p}=1$.
\end{itemize}
The quadratic term of \eqref{bck} is $zF_p^2$.
\epf

\subsection{An a priori bound}
The following elementary lemma ensures that the Borel transforms we consider are entire.

\bl\label{lem:entire}
Let $0<\mu<1$ and $S\ge1$, and let $\pi_1,\dots,\pi_S$ be polynomials with
$\deg\pi_s\le\mu s$. If $x_n=\sum_{s=1}^S\pi_s(n)\,x_{n-s}$ for all $n\ge N$, then there are
constants $M,K$ such that $|x_n|\le MK^n(n!)^{\mu}$ for all $n\ge0$. In particular
$\sum_nx_nz^n/n!$ is an entire function.
\el

\bpf
There is $B$ with $|\pi_s(n)|\le Bn^{\mu s}$ for $n\ge1$ and $1\le s\le S$. For $n\ge2S$ we have
$(n)_s\ge(n/2)^s$, where $(n)_s=n(n-1)\cdots(n-s+1)$. Hence
$|\pi_s(n)|K^{-s}(n)_s^{-\mu}\le B(2/K)^s$. Choose $K\ge1$ with $\sum_{s=1}^SB(2/K)^s\le1$, let
$N'=\max(N,2S)$ and $M=\max_{n<N'}|x_n|/(K^n(n!)^\mu)$. Dividing the recurrence by
$K^n(n!)^\mu$, induction on $n\ge N'$ gives $|x_n|/(K^n(n!)^\mu)\le M$. Finally,
$|x_n|/n!\le MK^n(n!)^{\mu-1}$, and $(n!)^{1-\mu}$ grows faster than any exponential.
\epf

\subsection{Linear differential equations with an irregular singularity at infinity}
\label{sec:HT}
We collect the classical facts that we use; see \cite{Wa87} and \cite{Ba00}. Let
$L=\sum_{i=0}^NP_i(z)D^i$, where $D=\mathrm d/\mathrm dz$, be a differential operator with
polynomial coefficients and constant $P_N\ne0$. Represent each monomial $z^dD^i$ occurring in $L$
by the point $(i,d)$. The \emph{Newton polygon} of $L$ is the upper boundary of the convex hull of
these points, starting at $(N,0)$. The slope of an edge joining $(i_1,d_1)$ to $(i_2,d_2)$, with
$i_1>i_2$, is $\gamma=(d_2-d_1)/(i_1-i_2)$, and its length is $i_1-i_2$.

\begin{itemize}
\item[(F1)] (Hukuhara--Turrittin) $Ly=0$ has a formal fundamental system of $N$ formal solutions
\[
\hat y_k=e^{Q_k(z)}z^{\rho_k}\sum_{j=0}^{J_k}(\log z)^j\,\hat s_{k,j}(z^{-1/s}),
\]
where $s\ge1$ is an integer, $Q_k$ is a polynomial in $z^{1/s}$ without constant term,
$\rho_k\in\mathbb C$, and the $\hat s_{k,j}$ are formal power series. If $\gamma_{\max}>0$ is the
largest slope of the Newton polygon and $\ell$ the length of the corresponding edge, then exactly
$\ell$ of the $Q_k$ have degree $1+\gamma_{\max}$ in $z$. All other $Q_k$ vanish or have degree
$<1+\gamma_{\max}$.
\item[(F2)] (Existence in sectors) For every direction $\theta_0$ there are $\delta>0$ and $R>0$
such that, in $S(\theta_0)=\{|\arg z-\theta_0|<\delta,\ |z|>R\}$, there are linearly independent
solutions $y_1,\dots,y_N$ of $Ly=0$ with $y_k\sim\hat y_k$. This means that
$e^{-Q_k}z^{-\rho_k}y_k$ has the asymptotic expansion $\sum_j(\log z)^j\hat s_{k,j}(z^{-1/s})$,
uniformly in $S(\theta_0)$, for a continuous choice of the branches of $z^{1/s}$ and $\log z$ in
$S(\theta_0)$.
\end{itemize}

No condition on Stokes directions is needed in (F2): the sector $S(\theta_0)$ may contain
directions along which two exponential parts have equal real parts. What fails there is the
uniqueness of the $y_k$, and the bases attached to overlapping sectors differ by constant matrices
(Stokes phenomenon). We only use the existence statement~(F2), never uniqueness or connection
formulas. It has the following consequence.

\begin{itemize}
\item[(F3)] Let $y$ be a solution of $Ly=0$ holomorphic in $\{|z|>R'\}$, and let $\Theta$ be a
compact set of directions. Then there are $M,K,R$ such that
$|y(re^{i\theta})|\le Mr^K\exp\bigl(\max_k\operatorname{Re}Q_k(re^{i\theta})\bigr)$ for
$\theta\in\Theta$ and $r\ge R$. For each $\theta$, the maximum is taken over the exponential parts
with the branches attached to a sector containing $\theta$.
\end{itemize}

Indeed, $\Theta$ is covered by finitely many sectors as in (F2). In each of them, $y$ is a
constant linear combination of the corresponding basis, and each basis element satisfies such a
bound.

\section{The linearized sequence for $p=2$}\label{sec:lin2}

In Sections~\ref{sec:lin2} to~\ref{sec:f2} we set $p=2$, $G=G_2$ and $g_n=g^{(2)}_n$. Thus $G$ is
the solution of
\begin{equation}\label{linearizedeq}
G(z)=z[u^2]\frac{M(z,u)}{1-u} + z [u^2]\frac1{1-u}G\left(\frac{z}{1-2zM(z,u)}\right).
\end{equation}

\bl\label{lem:rec}
The coefficients $g_n=[z^n]G(z)$ satisfy the recurrence relation
\begin{multline}\label{rec}
g_n=4g_{n-1}-6g_{n-2}+(2n-2)g_{n-3}-(4n-15)g_{n-4}+(2n^2-14n+20)g_{n-5} \\
-(2n^2-22n+60)g_{n-6}, \quad\mbox{ for } n\ge 6,
\end{multline}
with initial condition $g_0=g_1=0$, $g_2=1$, $g_3=2$, $g_4=4$, $g_5=12$.
\el

\bpf
If we extract the coefficients in Equation~\eqref{linearizedeq} and use that the partial
derivatives of $M(z,u)$ satisfy $M(z,0)=0$, $M_u(z,0)=z/(1-z)$ and $M_{uu}(z,0)=2z^3/(1-z)^3$,
then we obtain
\[
G(z)=\frac{z^2(2z^2-2z+1)}{(1-z)^3}+zG(z)+\frac{2z^4(1-2z+4z^2-2z^3)}{(1-z)^3} G'(z)
+\frac{2z^7}{(1-z)^2} G''(z).
\]
Multiplying by $(1-z)^3$ gives
\[
(1-z)^3G-z(1-z)^3G-2z^4(1-2z+4z^2-2z^3)G'-2z^7(1-z)G''=z^2(1-2z+2z^2).
\]
We extract the coefficient of $z^n$, using $[z^n]z^aG^{(k)}=(n-a+1)\cdots(n-a+k)\,g_{n-a+k}$. For
$n\ge5$ the right-hand side does not contribute. Collecting the terms, the coefficients of
$g_n,g_{n-1},\dots,g_{n-6}$ are $1$, $-4$, $6$, $-2(n-1)$, $4n-15$, $-2(n-2)(n-5)$ and
$2(n-5)(n-6)$, which is \eqref{rec}; it holds for $n\ge5$. The equation above is the case $p=2$ of Lemma~\ref{lem:linp}.
\epf

\brem
Let us also remark that generating the first few terms and using guess-and-prove by means of gfun
we can obtain several recurrence relations for the coefficients of $G(z)$. Indeed, this depends
on the preassigned parameters of the guess-and-prove procedure (which are the bounds for the order
of the recurrence itself as well as the degree of its coefficients). If we aim at finding a
smallest order recurrence, we get, for $n\ge0$,
\begin{multline*}
(n+1)\, g_{n+5}=2n(n+1)(n+2)\,g_n-2(n+1)(n+2)\,g_{n+1}+(n+2)(2n+3)\,g_{n+2}
\\ -(3n+5)\,g_{n+3}+(3n+4)\,g_{n+4}.
\end{multline*}
But for the arguments we use in the following, the recurrence given in Equation~\eqref{rec} was
more convenient.
\erem

\brem\label{rem:frec}
For comparison, Proposition~\ref{BCKdirect} yields for $f_n$ a recurrence with a quadratic
term and unbounded order. The sequence $f_n$ is the solution of the recurrence relation
\begin{align*}
f_n&=m_n+\sum_{k=2}^{n-3} f_kf_{n-k-1}+f_{n-1}+(2n-6)f_{n-3}
+\sum_{k=4}^{n-2} \delta_{n,k}\, f_{n-k}, \\
\delta_{n,k}&=(2k-8)n^2-(3k^2-11k-6)n+k(k^2-3k-6),
\end{align*}
where $m_n$ is itself given by the recurrence $(n-2)^2m_n=(9n-26)m_{n-1}+(n-3)(n-6)m_{n-2}$
(equivalently, $m_n=1+\binom{n-2}{2}$ for $n\ge2$), and the initial conditions of these two
recurrences are $m_1= 0$, $m_2=1$ and $f_0=0$, $f_1=0$, $f_2=1$, $f_3=2$, $f_4=4$, $f_5=13$.

The polynomial $\delta_{n,k}$ can be obtained from the following formula given in \cite{BGGJ13}
and using creative telescoping methods \cite{Ze91}: for $4\le k\le n-2$ we have
$\delta_{n,k}=S(n,n-k)$, where
\[
S(n,l)=
\begin{cases}\displaystyle
\sum _{t=0}^{ \min\left(2,\left\lfloor \frac{n-l-1}2\right\rfloor\right)} \sum_{r=0}^t {\frac
{r{2}^{r}{l-1+r \choose r} \left( n-l-2-r \right) !}
{t!\, \left( t-r \right) !\, \left( n-l-1-2\,t \right) !}} & \mbox{ if } 1\le l<n-1,
\\[2mm]
1& \mbox{ if } l=n-1.
\end{cases}
\]
(For $k=1,2,3$ the same expression $S(n,n-k)$ gives $1$, $0$ and $2n-6$, i.e. the coefficients of
$f_{n-1}$, $f_{n-2}$ and $f_{n-3}$ above.)
\erem

Intuitively, the recurrence tells us that the fifth term on the right-hand side is dominant and
thus the growth rate should be of order $(n!)^{2/5}$. The following lemmas make this precise.

\bl\label{ratio2}
For all $n\ge 0$ we have $g_{n+1}\ge 2g_n$. In particular $g_{n+1}/g_n\ge 2$ for $n\ge 2$.
\el

\bpf
We set $g_k=0$ for $k<0$. The functional equation used in the proof of Lemma~\ref{lem:rec} reads
\[
G(z)=A(z)+zG(z)+B_1(z)G'(z)+B_2(z)G''(z)
\]
with
\[
A(z)=\frac{z^2}{1-z}+\frac{z^4}{(1-z)^3},\quad
B_1(z)=\frac{2z^4}{1-z}+\frac{2z^6(3-2z)}{(1-z)^3},\quad
B_2(z)=\frac{2z^7}{(1-z)^2}.
\]
(Indeed $A(z)=z[u^2]M(z,u)/(1-u)$, and $B_1(z)$ coincides with the coefficient of $G'(z)$ given in
the proof of Lemma~\ref{lem:rec}.) Write $A(z)=\sum_n a_nz^n$, $B_1(z)=\sum_j\beta_jz^j$ and
$B_2(z)=\sum_j\gamma_jz^j$. Then $a_0=a_1=0$ and $a_n=1+\binom{n-2}{2}$ for $n\ge2$. Moreover
$\beta_j\ge0$ for all $j$, since $(3-2z)/(1-z)^3=\sum_{k\ge0}\frac{(k+1)(k+6)}2z^k$; also
$\beta_j=0$ for $j<4$ and $\beta_4=2$. Finally $\gamma_j=2(j-6)\ge0$ for $j\ge7$ and $\gamma_j=0$
otherwise. Extracting coefficients gives, for all $n\ge0$,
\begin{equation}\label{posform}
g_n=a_n+g_{n-1}+\sum_{j\ge4}\beta_j\,(n-j+1)\,g_{n-j+1}
+\sum_{j\ge7}\gamma_j\,(n-j+2)(n-j+1)\,g_{n-j+2},
\end{equation}
where all coefficients are nonnegative.

We prove $g_{k+1}\ge 2g_k$ by induction on $k$. It holds trivially for $k<0$, and the initial
values $g_0,\dots,g_4=0,0,1,2,4$ give it for $0\le k\le3$. Let $n\ge4$ and assume
$g_{k+1}\ge2g_k$ for all $k<n$. Subtracting twice \eqref{posform} at $n$ from \eqref{posform} at
$n+1$ yields
\[
g_{n+1}-2g_n=(a_{n+1}-2a_n)+(g_n-2g_{n-1})+\sum_{j\ge4}\beta_jD_j+\sum_{j\ge7}\gamma_jE_j,
\]
where, with $m=n-j+1$ in $D_j$ and $m=n-j+2$ in $E_j$,
\[
D_j=(m+1)g_{m+1}-2m\,g_m,\qquad E_j=m\bigl((m+1)g_{m+1}-2(m-1)g_m\bigr).
\]
In both cases $m\le n-3<n$, so $2g_m\le g_{m+1}$ by the induction hypothesis. For $m\ge0$ this
gives $D_j\ge (m+1)g_{m+1}-m\,g_{m+1}=g_{m+1}\ge0$ and $E_j\ge 2m\,g_{m+1}\ge0$. For $m<0$ we
have $g_m=0$ and $(m+1)g_{m+1}=0$, so $D_j=E_j=0$. All terms of both sums are therefore
nonnegative; keeping only the term $j=4$ of the first sum, for which $\beta_4=2$ and
$D_4\ge g_{n-2}$, we get
\[
g_{n+1}-2g_n\ge (a_{n+1}-2a_n)+(g_n-2g_{n-1})+2g_{n-2}.
\]
By the induction hypothesis $g_n-2g_{n-1}\ge0$, and $g_{n-2}\ge 2^{n-4}g_2=2^{n-4}$. A direct
computation gives $a_{n+1}-2a_n=-1-\frac{(n-2)(n-5)}2$ for $n\ge2$. Hence
\[
g_{n+1}-2g_n\ge 2^{n-3}-1-\frac{(n-2)(n-5)}2.
\]
The right-hand side equals $2$ for $n=4$, and it does not decrease when $n$ increases by one,
because $2^{n-3}$ increases by $2^{n-3}$ while $\frac{(n-2)(n-5)}2$ increases by
$n-3\le2^{n-3}$. Therefore $g_{n+1}\ge2g_n$, which completes the induction.
\epf

\bcor\label{cor:n2}
For $n\ge 8$ we have $g_n\ge (n-2)(n-5)\,g_{n-5}$.
\ecor

\bpf
The result is immediate from Lemma~\ref{ratio2}. Indeed, for $n\ge8$ we have
\begin{align*}
&4g_{n-1}-6g_{n-2}\ge 2g_{n-2}\ge 0, \\
&(2n-2)g_{n-3}-(4n-15)g_{n-4}\ge 11g_{n-4} \ge 0, \\
&(2n^2-14n+20)g_{n-5} -(2n^2-22n+60)g_{n-6}\\
&\qquad\ge (n-2)(n-5)g_{n-5}+8(n-5)g_{n-6}\ge (n-2)(n-5)g_{n-5},
\end{align*}
where in the last line we wrote $2n^2-14n+20=2(n-2)(n-5)$, used
$(n-2)(n-5)g_{n-5}\ge 2(n-2)(n-5)g_{n-6}$ and $2(n-2)(n-5)-2(n-5)(n-6)=8(n-5)$.
\epf

\bcor\label{lower}
There is a constant $c>0$ such that $g_n\ge c\, n^{-7/5}(n!)^{2/5}$ for all $n\ge 3$.
\ecor

\bpf
Let $m_j=n-5j$ and let $J$ be the smallest index with $m_J\le 7$, so that $3\le m_J\le7$ and
$m_j\ge 8$ for $j<J$. Iterating Corollary~\ref{cor:n2} gives
$g_n\ge g_{m_J}\prod_{j<J}(m_j-2)(m_j-5)$, while
\[
n!=m_J!\prod_{j<J}m_j(m_j-1)\cdots(m_j-4)\le 7!\prod_{j<J}m_j^5 .
\]
Since $g_{m_J}\ge2$,
\[
\frac{g_n}{(n!)^{2/5}}\ge \frac{2}{(7!)^{2/5}}\prod_{j<J}\left(1-\frac2{m_j}\right)\left(1-\frac5{m_j}\right).
\]
Using $\log(1-x)\ge -x-x^2$ for $0\le x\le 5/8$, the product is at least
\[
\exp\Bigl(-7\sum_{j<J}m_j^{-1}-29\sum_{j<J}m_j^{-2}\Bigr).
\]
Finally
$\sum_{j<J}m_j^{-1}\le\frac15\int_{3}^{n}\frac{dx}{x}\le\frac15\log n$ and
$\sum_{j<J}m_j^{-2}\le\sum_{m\ge8}m^{-2}<1/7$, which proves the claim.
\epf

Thus we are faced with a rapidly growing sequence, as expected.

\bl\label{upper}
For every $\eps>0$ there is a constant $C_\eps$ such that the sequence $(g_n)_{n\ge 0}$ satisfies
$g_n\le C_\eps\,(n!)^{2/5}(2+\eps)^{n/5}$ for all $n\ge0$.
\el

\bpf
We choose an arbitrary $\eps>0$ and substitute
\[
v_n=v_n(\eps):=\frac{g_n}{(n!)^{2/5} (2+\eps)^{n/5}}.
\]
Then Equation~\eqref{rec} becomes, with $(n)_k=n(n-1)\cdots(n-k+1)$,
\begin{align*}
v_n &=\sum_{k=1}^6 a_k(n) v_{n-k} \\
&=\frac{4}{n^{2/5} (2+\eps)^{1/5}} v_{n-1} - \frac{6}{(n)_2^{2/5}(2+\eps)^{2/5}} v_{n-2}
+\frac{2n-2}{(n)_3^{2/5}(2+\eps)^{3/5}} v_{n-3} \\
&\qquad -\frac{4n-15}{(n)_4^{2/5}(2+\eps)^{4/5}} v_{n-4}
+\frac{2n^2-14n+20}{(n)_5^{2/5}(2+\eps)} v_{n-5}
-\frac{2n^2-22n+60}{ (n)_6^{2/5}(2+\eps)^{6/5}} v_{n-6}.
\end{align*}
For $n\ge7$ the coefficients $a_2,a_4,a_6$ are nonpositive and $v_j\ge0$, hence
$v_n\le (a_1(n)+a_3(n)+a_5(n))\max_{n-5\le j\le n-1}v_j$. As $n\to\infty$, $a_1(n)\to0$,
$a_3(n)\to0$ and $a_5(n)\to 2/(2+\eps)<1$, so there is $n_0(\eps)$ with
$a_1(n)+a_3(n)+a_5(n)\le1$ for $n\ge n_0(\eps)$. By induction, $v_n\le
C_\eps:=\max_{j<n_0(\eps)}v_j$ for all $n$.
\epf

\begin{figure}[htb]
\begin{center}
\includegraphics[width=0.55\textwidth]{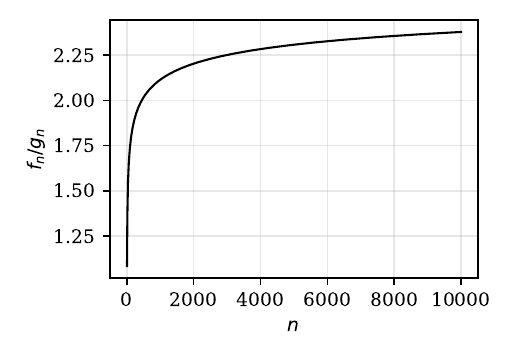}
\end{center}
\caption{\small The quotient $f_n/g_n$ for $n\le10000$.}
\label{fig1}
\end{figure}

In order to get information on the sequence $(f_n)_{n\ge 0}$ one might be tempted to add the
quadratic term to the recurrence, but this is not feasible. But it may serve as a heuristic.
Indeed putting the quadratic term and substituting as before, this term becomes
$(2+\eps)^{-1/5} \sum_{l=1}^{n-1}\bigl(n{n-1\choose l}\bigr)^{-2/5} v_l v_{n-1-l}$ which is
obviously very small. This indicates that the quadratic term is not essential for the asymptotic
behaviour of $f_n$. We will show that it actually only affects the multiplicative factor of the
asymptotic expression. The sequences start with
\begin{align*}
(f_n)_{n\ge 0}&=(0,0,1,2,4,13,40,127,459,1671,6284,\dots ),\\
(g_n)_{n\ge 0}&=(0,0,1,2,4,12,35,110,390,1384,5101,\dots ),
\end{align*}
and Figure~\ref{fig1} depicts the behaviour of the quotient for the first 10000 terms.

\section{The Borel transform for $p=2$}\label{sec:borel2}

Since we are facing a holonomic sequence we can use again the rich toolbox available for such
sequences. Due to the rapid growth of $g_n$, $G(z)$ is not an analytic function. Thus we will work
with its Borel transform
\[
H(z)=\sum_{n\ge 0} \frac{g_n}{n!} z^n
\]
which is an entire function. It satisfies the following differential equation.

\bl\label{lem:diffeq}
The function $y=H'(z)$ is a solution of the differential equation
\begin{equation}\label{diffeq}
y^{(iv)}=4y'''+(2z-6)y''+(2z^2-4z+8)y'+(-2z^2+8z-5)y
\end{equation}
with initial conditions $H(0)=H'(0)=0$, $H''(0)=1$, $H'''(0)=2$, $H^{(iv)}(0)=4$.
\el

\bpf
For a sequence $(x_n)$, let $R_n(x)$ denote the difference of the two sides of \eqref{rec} at
index $n$, with $x_k=0$ for $k<0$. Let $X(z)=\sum_nx_nz^n/n!$ and let $L_4$ be the operator such
that \eqref{diffeq} reads $L_4y=0$. A direct computation of coefficients shows
\[
L_4(X')=\sum_{k\ge0}R_{k+5}(x)\,\frac{z^k}{k!}.
\]
By the proof of Lemma~\ref{lem:rec}, $R_n(g)=0$ for $n\ge5$, hence $L_4(H')=0$. The initial
conditions are $g_0,\dots,g_4$.
\epf

\bp\label{fracpower}
Let
\[
Q(z)=\frac{3\cdot2^{1/3}}5z^{5/3}+2^{-4/3}z^{4/3}+z-\frac{2^{1/3}}{27}z^{2/3}+\frac{2^{5/3}}{81}z^{1/3},
\]
where principal branches are used. There exist $\delta>0$ and a constant $C>0$ such that, as
$z\to\infty$ in the sector
\begin{equation}\label{sector}
S=\{z\in\mathbb C\;\mid\; |\arg(z)|<\delta\},
\end{equation}
we have, uniformly in $\arg z$,
\begin{equation}\label{H-asym}
H(z)=\frac{C}{z}\,e^{Q(z)}\left(1+O\bigl(|z|^{-1/3}\bigr)\right).
\end{equation}
\ep

\bpf
\emph{The equations.} Denote by $L_4$ the operator of \eqref{diffeq}, i.e.
$L_4y=y^{(iv)}-4y'''-(2z-6)y''-(2z^2-4z+8)y'-(-2z^2+8z-5)y$, and set $L_5h=L_4(h')$. By
Lemma~\ref{lem:diffeq}, $H$ satisfies the fifth order equation $L_5H=0$. On any simply connected
domain the solutions of $L_5h=0$ form a vector space of dimension $5$ which contains the
constants.

\emph{Formal solutions of $L_4y=0$.} Writing $y=\exp\left(\int w\right)$ we have
$L_4y=y\,\mathcal R(w)$ with
\begin{multline*}
\mathcal R(w)=w'''+4ww''+3w'^2+6w^2w'+w^4-4\left(w''+3ww'+w^3\right)\\
-(2z-6)(w'+w^2)-(2z^2-4z+8)w+2z^2-8z+5 .
\end{multline*}
Let $w_0$ be one of the three cube roots of $2$ and look for a formal series
$w=\sum_{j\ge0}w_jz^{(2-j)/3}$. The coefficient of $z^{8/3}$ in $\mathcal R(w)$ is
$w_0(w_0^3-2)=0$. For $j\ge1$, the unknown $w_j$ appears in the coefficient of $z^{(8-j)/3}$ only
through $w^4$ and $-2z^2w$ (all other terms contain $w_j$ multiplied by a lower power of $z$), so
this coefficient equals $(4w_0^3-2)w_j+P_j=6w_j+P_j$, where $P_j$ is a polynomial in
$w_0,\dots,w_{j-1}$. Hence the $w_j$ are uniquely determined recursively. Integrating term by
term, the term $j=5$ produces a logarithm and the terms $j\ge6$ produce negative powers, so
\[
\hat y=\exp\left(\int w\right)=e^{P(z)}z^{w_5}\left(1+\sum_{k\ge1}s_kz^{-k/3}\right),\qquad
P(z)=\sum_{j=0}^4\frac{3w_j}{5-j}z^{(5-j)/3}.
\]
For $w_0=2^{1/3}$ the recursion gives $w_1=2^{2/3}/3$, $w_2=1$, $w_3=-2^{4/3}/81$,
$w_4=2^{5/3}/243$ and $w_5=-1/3$, so $P=Q$ and we obtain a formal solution
$\hat y_1=e^{Q(z)}z^{-1/3}(1+O(z^{-1/3}))$. The equation $\mathcal R(w)=0$ only involves $z=t^3$
and $\mathrm d/\mathrm dz$, which are invariant under $t\mapsto \omega t$, $\omega=e^{2\pi i/3}$.
Hence the two other choices of $w_0$ yield the formal solutions $\hat y_2,\hat y_3$ obtained from
$\hat y_1$ by replacing $z^{1/3}$ with $\omega^{\pm1}z^{1/3}$; their exponential parts $Q_2,Q_3$
have leading terms $\frac{3\cdot2^{1/3}}5\omega^{\mp1}z^{5/3}$.

A fourth formal solution is obtained by writing $y=e^zu$: the equation becomes
$u^{(iv)}-2zu''-(2z^2+4)u'-6zu=0$. Inserting $u=z^{-3}\bigl(1+\sum_{k\ge1}c_kz^{-k}\bigr)$, the
coefficient of $c_k$ at the leading order $z^{-2-k}$ is $-2(-3-k)-6=2k\neq0$, while all the other
terms are of lower order. So the $c_k$ are determined recursively and
$\hat y_4=e^{z}z^{-3}(1+O(z^{-1}))$.

\emph{Formal solutions of $L_5h=0$.} Let $\hat y=e^{P}z^{\rho}\hat s$ be one of the
$\hat y_k$, with $\hat s=1+O(z^{-1/3})$ and $P'(z)=pz^{\mu}(1+O(z^{-1/3}))$, $p\ne0$, $\mu\ge0$.
We look for $\hat h=e^{P}z^{\rho-\mu}\hat v$ with $\hat h'=\hat y$, i.e.
\[
\left(z^{-\mu}P'+(\rho-\mu)z^{-\mu-1}\right)\hat v+z^{-\mu}\hat v'=\hat s .
\]
The factor $z^{-\mu}P'+(\rho-\mu)z^{-\mu-1}=p(1+O(z^{-1/3}))$ is invertible, and
$\hat v\mapsto z^{-\mu}\hat v'$ raises the valuation in $z^{-1/3}$ by at least $3$. Hence this
equation has a unique solution $\hat v=1/p+O(z^{-1/3})$, obtained as the limit of the iteration
$\hat v\mapsto\bigl(z^{-\mu}P'+(\rho-\mu)z^{-\mu-1}\bigr)^{-1}(\hat s-z^{-\mu}\hat v')$ in the
$z^{-1/3}$-adic topology. Since $L_5\hat h=L_4\hat y=0$, $\hat h$ is a formal solution of
$L_5h=0$. For $\hat y_1$ we have $P=Q$, $\mu=2/3$, $p=2^{1/3}$, $\rho=-1/3$, thus
\[
\hat h_1=2^{-1/3}z^{-1}e^{Q(z)}\left(1+O(z^{-1/3})\right),
\]
and similarly we obtain $\hat h_2,\hat h_3$ (exponential parts $Q_2,Q_3$) and
$\hat h_4=e^{z}z^{-3}(1+O(z^{-1/3}))$. Together with $\hat h_0=1$, these are five formal solutions
with pairwise distinct exponential parts $0,Q,Q_2,Q_3,z$; they are therefore linearly independent
and form a formal fundamental system of $L_5h=0$.

\emph{Actual solutions.} By (F2) of Section~\ref{sec:HT}, applied with $\theta_0=0$, there are $\delta_0>0$, $R>0$ and solutions $h_1,\dots,h_4$ of $L_5h=0$, holomorphic in
$S'=\{|\arg z|<\delta_0,\ |z|>R\}$, such that $h_k\sim\hat h_k$ uniformly there; in particular
\[
h_1(z)=2^{-1/3}z^{-1}e^{Q(z)}\left(1+O(|z|^{-1/3})\right),
\]
$|h_{2}(z)|,|h_3(z)|\le |z|^{K}\exp(\operatorname{Re}Q_{2,3}(z))$ and
$|h_4(z)|\le|z|^{K}e^{\operatorname{Re}z}$ for some $K$. Since $H$ is entire and $L_5H=0$, there
are constants $c_0,\dots,c_4$ with $H=c_0+c_1h_1+c_2h_2+c_3h_3+c_4h_4$ in $S'$.

\emph{Domination.} Let $\delta=\min(\delta_0,\pi/10)$ and $z=re^{i\theta}$ with
$|\theta|<\delta$. Then $|5\theta/3|<\pi/6$, hence $\cos(5\theta/3)>\sqrt3/2$ while
$\cos(5\theta/3\mp2\pi/3)\le0$. Therefore, with $a=\frac{3\cdot2^{1/3}}5$,
\[
\operatorname{Re}Q(z)\ge \tfrac{\sqrt3}{2}\,a\,r^{5/3}-O(r^{4/3}),\qquad
\operatorname{Re}Q_{2,3}(z)\le O(r^{4/3}),\qquad \operatorname{Re}z\le r,
\]
and consequently $|h_k(z)/h_1(z)|$ and $|1/h_1(z)|$ are $O\bigl(\exp(-\tfrac{\sqrt3}{4}ar^{5/3})\bigr)$
uniformly in $S$ for $k=2,3,4$. If $c_1\ne0$, this gives
$H(z)=c_1h_1(z)\bigl(1+O(e^{-\frac{\sqrt3}{4}ar^{5/3}})\bigr)$, i.e. \eqref{H-asym} with
$C=2^{-1/3}c_1$.

\emph{Non-vanishing of $c_1$ and positivity of $C$.} Suppose $c_1=0$. On the positive real axis, $h_2,h_3$ are
bounded and $h_4(r)=O(r^Ke^r)$, so $H(r)=O(r^Ke^r)$. But Lemma~\ref{ratio2} gives
$g_n\ge2^{n-2}$ for $n\ge2$, hence
\[
H(r)\ge\sum_{n\ge2}2^{n-2}\frac{r^n}{n!}=\frac{e^{2r}-1-2r}4,
\]
a contradiction. Thus $c_1\neq0$. Finally, by the domination estimate,
\[
C=2^{-1/3}c_1=\lim_{r\to+\infty}rH(r)e^{-Q(r)},
\]
which is a limit of positive real numbers since $H$ has nonnegative coefficients and
$Q(r)\in\mathbb R$. Hence $C\ge0$, and $C\ne0$, so $C>0$.
\epf

Knowing the asymptotic behaviour of the function, the next step is to derive the asymptotic
behaviour of its coefficients. Due to the rapid growth of the function \cite[Sec.~12]{Od95}
suggests the application of the saddle point method. Indeed the function fits well into the
framework of Hayman admissible functions \cite{Hayman}.

\bl\label{lem:hayman}
The function $H(z)$ is Hayman-admissible. More precisely, let $a(r)=rH'(r)/H(r)$,
$b(r)=ra'(r)$ and $\eta(r)=r^{-2/3}$. Then, as $r\to\infty$,
\begin{enumerate}
\item[(i)] $b(r)\to\infty$;
\item[(ii)] $H(re^{i\phi})=H(r)\,e^{i\phi a(r)-\phi^2b(r)/2}\,(1+o(1))$ uniformly for
$|\phi|\le\eta(r)$;
\item[(iii)] $H(re^{i\phi})=o\bigl(H(r)/\sqrt{b(r)}\bigr)$ uniformly for
$\eta(r)\le|\phi|\le\pi$.
\end{enumerate}
\el

\bpf
Write $Q(z)=\sum_{j=1}^5\alpha_jz^{j/3}$ with $\alpha_5=\frac{3\cdot2^{1/3}}5$,
$\alpha_4=2^{-4/3}$, $\alpha_3=1$, $\alpha_2=-\frac{2^{1/3}}{27}$, $\alpha_1=\frac{2^{5/3}}{81}$.
Note that $\alpha_j>0$ for $j\ne2$. Let $\delta$, $C$ and $S$ be as in
Proposition~\ref{fracpower}; decreasing $\delta$ if necessary, we may assume $\delta\le\pi/10$.

\emph{(a) Logarithmic derivatives.} By \eqref{H-asym}, $H$ does not vanish in $S$ for $|z|$ large,
and $\Lambda(z)=\log H(z)=\log C-\log z+Q(z)+E(z)$, where $E$ is holomorphic and
$E(z)=O(|z|^{-1/3})$ uniformly in $S$. The disk of radius $|z|\sin(\delta/2)$ centred at a point
$z$ with $|\arg z|\le\delta/2$ lies in $S$, so Cauchy's estimates give
$E^{(k)}(z)=O(|z|^{-1/3-k})$ for $k=1,2,3$, uniformly for $|\arg z|\le\delta/2$. Since
$Q'(z)\asymp|z|^{2/3}$, $Q''(z)\asymp|z|^{-1/3}$ and $Q'''(z)\asymp|z|^{-4/3}$, we obtain
\begin{align*}
a(r)&=rQ'(r)+O(1)=2^{1/3}r^{5/3}\bigl(1+O(r^{-1/3})\bigr),\\
b(r)&=r\bigl(rQ'(r)\bigr)'+O(r^{-1/3})=\tfrac53\,2^{1/3}r^{5/3}\bigl(1+O(r^{-1/3})\bigr).
\end{align*}
This proves (i).

\emph{(b) Proof of (ii).} Let $f(\phi)=\Lambda(re^{i\phi})$. Since $\mathrm d/\mathrm d\phi=
iz\,\mathrm d/\mathrm dz$, we have $f(0)=\log H(r)$, $f'(0)=ia(r)$, $f''(0)=-b(r)$ and
$f'''(\phi)=-i\bigl(z\Lambda'(z)+3z^2\Lambda''(z)+z^3\Lambda'''(z)\bigr)=O(r^{5/3})$ uniformly for
$|\phi|\le\delta/2$, by (a). Taylor's formula gives
$f(\phi)=\log H(r)+i\phi a(r)-\phi^2b(r)/2+O(r^{5/3}|\phi|^3)$, and the error term is
$O(r^{-1/3})$ for $|\phi|\le r^{-2/3}$.

\emph{(c) Proof of (iii) for $r^{-2/3}\le|\phi|\le\delta$.} By \eqref{H-asym},
\[
\frac{|H(re^{i\phi})|}{H(r)}=\exp\bigl(\operatorname{Re}Q(re^{i\phi})-Q(r)\bigr)\bigl(1+O(r^{-1/3})\bigr).
\]
Using
$\alpha_j>0$ for $j\ne2$ and $\frac{2x^2}{\pi^2}\le1-\cos x\le\frac{x^2}2$ for $|x|\le\pi$,
\begin{align*}
\operatorname{Re}Q(re^{i\phi})-Q(r)&=-\sum_{j=1}^5\alpha_jr^{j/3}\Bigl(1-\cos\frac{j\phi}3\Bigr)\\
&\le-\frac{50\,\alpha_5}{9\pi^2}r^{5/3}\phi^2+\frac{2|\alpha_2|}{9}r^{2/3}\phi^2\le-c_1r^{5/3}\phi^2
\end{align*}
for some $c_1>0$ and $r$ large. For $|\phi|\ge r^{-2/3}$ this is at most $-c_1r^{1/3}$, hence
$|H(re^{i\phi})|\le2H(r)e^{-c_1r^{1/3}}=o\bigl(H(r)/\sqrt{b(r)}\bigr)$, since $b(r)\asymp r^{5/3}$.

\emph{(d) Proof of (iii) for $\delta\le|\phi|\le\pi$.} Recall from the proof of
Proposition~\ref{fracpower} that $L_5h=0$ has a formal fundamental system with exponential parts
$0$, $z$ and $Q_1=Q,Q_2,Q_3$, the latter being obtained from $Q$ by taking the three
determinations of $z^{1/3}$. For $z=re^{i\theta}$ the set
$\{\operatorname{Re}Q_k(z):k=1,2,3\}$ is therefore
\[
\{F_r(\psi):\psi\in\Psi(\theta)\},\quad\text{where}\quad F_r(\psi)=\sum_{j=1}^5\alpha_jr^{j/3}\cos(j\psi)
\]
and $\Psi(\theta)=\bigl\{\frac{\theta+2\pi m}{3}:m\in\{-1,0,1\}\bigr\}$,
independently of the choice of branch. Applying (F3) of Section~\ref{sec:HT} to $H$ and to the
arc $A=\{\delta\le|\theta|\le\pi\}$, we obtain constants $M,K,R$ such that
\begin{equation}\label{hayman_global}
|H(re^{i\theta})|\le M\,r^K\exp\Bigl(\max\Bigl\{0,\ r\cos\theta,\ \max_{\psi\in\Psi(\theta)}F_r(\psi)\Bigr\}\Bigr)
\end{equation}
for $r\ge R$ and $\theta\in A$.
For $\theta\in A$, every $\psi\in\Psi(\theta)$ satisfies $\delta/3\le|\psi|\le\pi$. The function
$\mu(\psi)=\max(1-\cos5\psi,\,1-\cos4\psi)$ vanishes only if $5\psi,4\psi\in2\pi\mathbb Z$, i.e.
$\psi\in2\pi\mathbb Z$. Hence $\mu_\delta=\min\{\mu(\psi):\delta/3\le|\psi|\le\pi\}>0$. Keeping
only the terms $j=4,5$ of the nonnegative terms, we get for $r\ge1$ and $\psi\in\Psi(\theta)$
\[
F_r(\psi)=Q(r)-\sum_{j=1}^5\alpha_jr^{j/3}(1-\cos j\psi)
\le Q(r)-\min(\alpha_4,\alpha_5)\,\mu_\delta\,r^{4/3}+2|\alpha_2|r^{2/3}.
\]
Since moreover $\max(0,r\cos\theta)\le r$, \eqref{hayman_global} yields
$|H(re^{i\theta})|\le Mr^K\exp\bigl(Q(r)-c_2r^{4/3}\bigr)$ for some $c_2>0$ and $r$ large. On the
other hand $H(r)\ge\frac C2r^{-1}e^{Q(r)}$ by \eqref{H-asym}. Therefore
\[
\frac{|H(re^{i\theta})|}{H(r)}\le\frac{2M}{C}\,r^{K+1}e^{-c_2r^{4/3}}=o\bigl(b(r)^{-1/2}\bigr)
\]
uniformly for $\theta\in A$. Together with (c), this proves (iii).
\epf

\brem
The term $r^{4/3}$ in (d) cannot be improved to $r^{5/3}$: for $\psi=\pm2\pi/5$ or $\pm4\pi/5$ we
have $\cos5\psi=1$, so the leading term alone does not separate these directions from the positive
axis; it is the second coefficient $\alpha_4>0$ that does.
\erem

By Hayman's theorem (see \cite{Hayman} or \cite[Sec.~VIII.5]{FlSe09}), Lemma~\ref{lem:hayman}
implies
\[
h_n=[z^n]H(z)\sim\frac{H(\zeta)}{\zeta^n\sqrt{2\pi b(\zeta)}},\qquad\text{where }a(\zeta)=n,
\]
and the positive solution $\zeta=z_n$ of $a(\zeta)=n$ is unique since $a'(r)=b(r)/r>0$.

\bp
The function $H(z)/z^n$ has a uniquely determined saddle point $z_n$ on the positive real axis
which satisfies
\[
z_n\sim 2^{-1/5}n^{3/5}-\frac{2^{1/5}}5 n^{2/5}-\frac{13\cdot2^{3/5}}{50}n^{1/5}+\frac15
+\frac{41\cdot 2^{2/5}}{625} n^{-1/5}+\frac{8069\cdot 2^{4/5}}{31250}n^{-2/5},
\]
as $n\to \infty$.
\ep

\bpf
The saddle point is the positive solution of the equation $zH'(z)/H(z)=n$. From \eqref{H-asym} we
can derive easily the assertion.
\epf

\bcor\label{cor:gasym}
Let
\begin{multline*}
\Phi(n)=n^{2n/5}2^{n/5}e^{-2n/5}n^{-3/5}\\
\times\exp\left(2^{-8/5}n^{4/5}+\frac{7\cdot 2^{4/5}}{15}n^{3/5}
-\frac{17\cdot 2^{1/5}}{75}n^{2/5}-\frac{41\cdot 2^{3/5}}{500}n^{1/5}\right).
\end{multline*}
There is a constant $A_g>0$ such that $g_n\sim A_g\Phi(n)$ as $n\to\infty$.
\ecor

\bpf
By Hayman's theorem and Lemma~\ref{lem:hayman},
\[
g_n=n!\,h_n\sim n!\,\frac{H(z_n)\,z_n^{-n}}{\sqrt{2\pi b(z_n)}}.
\]
We insert \eqref{H-asym},
$b(r)=\frac53 2^{1/3}r^{5/3}(1+O(r^{-1/3}))$, the expansion of $z_n$ and Stirling's formula. Since
$z_n$ is a stationary point of $\log H(r)-n\log r$, replacing it by its expansion truncated after
the term $n^{-2/5}$ (error $O(n^{-3/5})$) changes $\log(H(r)r^{-n})$ only by
$O\bigl(b(z_n)n^{-6/5}z_n^{-2}\bigr)=o(1)$. The remaining expansion is a routine computation
(carried out with computer algebra), which yields $g_n\sim A_g\Phi(n)$ with $A_g$ a positive
multiple of the constant $C$ of Proposition~\ref{fracpower}.
\epf

\section{From $g_n$ to $f_n$ for $p=2$}\label{sec:f2}

We now compare $f_n$ with $g_n$. The linear part of \eqref{bck} (for $p=2$) is the same as in
\eqref{linearizedeq}, so the computation in the proofs of Lemmas~\ref{lem:rec}
and~\ref{ratio2} gives $F=A+zF+B_1F'+B_2F''+zF^2$. For a sequence $x=(x_n)_{n\ge0}$ let
\[
(Tx)_n=x_{n-1}+\sum_{j\ge4}\beta_j\,(n-j+1)\,x_{n-j+1}+\sum_{j\ge7}\gamma_j\,(n-j+2)(n-j+1)\,x_{n-j+2}
\]
(with $x_k=0$ for $k<0$). Then \eqref{posform} reads $g=a+Tg$, and
\begin{equation}\label{posformF}
f=a+Tf+q,\qquad q_n=\sum_{i=2}^{n-3}f_if_{n-1-i}.
\end{equation}
The operator $T$ has nonnegative coefficients, and $(Tx)_n$ only depends on $x_k$ for $k<n$.
Hence, for every sequence $s$, the equation $x=s+Tx$ has a unique solution, which depends linearly
on $s$ and is nonnegative if $s$ is. In particular $f\ge g$, since $q\ge0$.

\bl\label{lem:green}
For $m\ge2$ let $e^{(m)}$ be the solution of $x=\delta_m+Tx$, where $\delta_m$ is the sequence
with a single $1$ at index $m$. Then:
\begin{enumerate}
\item[(i)] $f_n=g_n+\sum_{m=5}^nq_m\,e^{(m)}_n$ for all $n$;
\item[(ii)] $0\le e^{(m)}_n\le g_n/g_m$ for all $n$;
\item[(iii)] $e^{(m)}$ satisfies \eqref{rec} for all $n\ge m+4$.
\end{enumerate}
\el

\bpf
(i) follows from \eqref{posformF} by linearity, since $q_m=0$ for $m<5$ and $e^{(m)}_n=0$ for
$n<m$. (ii) We use induction on $n$. For $n<m$ we have $e^{(m)}_n=0$, and $e^{(m)}_m=1=g_m/g_m$.
For $n>m$, $e^{(m)}_n=(Te^{(m)})_n\le(Tg)_n/g_m\le(a+Tg)_n/g_m=g_n/g_m$, by the induction
hypothesis and the positivity of $T$. (iii) The generating function $E_m(z)=\sum_ne^{(m)}_nz^n$
satisfies $E_m=z^m+zE_m+B_1E_m'+B_2E_m''$. Multiplying by $(1-z)^3$ as in the proof of
Lemma~\ref{lem:rec}, the coefficient of $z^n$ gives exactly \eqref{rec} as soon as $n>m+3$, the
degree of $z^m(1-z)^3$.
\epf

\bl\label{lem:limit}
Let $(x_n)$ be a sequence satisfying \eqref{rec} for all $n\ge N$. Then $x_n/g_n$ converges to a
finite limit.
\el

\bpf
In \eqref{rec}, the coefficient of $x_{n-s}$ has degree $0,0,1,1,2,2$ for $s=1,\dots,6$, which
is at most $2s/5$. By Lemma~\ref{lem:entire}, $X(z)=\sum_nx_nz^n/n!$ is entire. Let $L_4$ and
$L_5$ be the operators of the proof of Proposition~\ref{fracpower}. By the identity in the proof
of Lemma~\ref{lem:diffeq}, $L_4(X')=\sum_{k\ge0}R_{k+5}(x)\,z^k/k!$, where only finitely many
$R_{k+5}(x)$ are nonzero. Hence $L_4(X')=P$ is a polynomial. Since $L_4(z^d)=2z^{d+2}+O(z^{d+1})$, we can
write $P=L_4(p_0)+c_0+c_1z$ with a polynomial $p_0$. Let $P_0'=p_0$ and $L_7=D^2\circ L_4\circ D$,
where $D=\mathrm d/\mathrm dz$. Then $L_7(X-P_0)=D^2(c_0+c_1z)=0$ and $L_7H=0$.

The equation $L_7h=0$ has the formal solutions $1,\hat h_1,\dots,\hat h_4$ constructed in the
proof of Proposition~\ref{fracpower}. It has two more: the formal antiderivatives of the formal series
$\hat u_0,\hat u_1$ in powers of $z^{-1}$ with $L_4\hat u_0=1$ and $L_4\hat u_1=z$. These series
exist because $L_4(z^{-k})=2z^{2-k}+O(z^{1-k})$, and their antiderivatives behave like
$-z^{-1}/2$ and $\frac12\log z$. These seven formal solutions are linearly independent, so they
form a formal fundamental system, with exponential parts $0$ (three times), $z$, $Q_1$, $Q_2$ and
$Q_3$. By (F2) (shrinking $\delta$ if necessary), there are solutions
$X_0,\dots,X_6$ of $L_7h=0$ in $S$ (for $|z|$ large) with $X_1\sim\hat h_1$, and
$|X_k(z)|\le M|z|^K\exp\bigl(\max(0,\operatorname{Re}z,\operatorname{Re}Q_2(z),\operatorname{Re}Q_3(z))\bigr)$
for $k\ne1$. Write $X-P_0=\sum_k\xi_kX_k$ and $H=\sum_k\eta_kX_k$ in $S$; by \eqref{H-asym},
$\eta_1\ne0$. Let $\lambda=\xi_1/\eta_1$ and $Y=X-\lambda H$.

In $S$, $Y=P_0+\sum_{k\ne1}(\xi_k-\lambda\eta_k)X_k$ has no component on $X_1$. By the domination estimate in
the proof of Proposition~\ref{fracpower}, the real parts of $Q_2$ and $Q_3$ are
$O(|z|^{4/3})$ in $S$, hence $|Y(z)|\le M|z|^K\exp(M|z|^{4/3})$ there. For
$\delta\le|\theta|\le\pi$, the argument of part~(d) of the proof of Lemma~\ref{lem:hayman} applies
verbatim to the solutions $X-P_0$ and $H$ of $L_7h=0$, because their exponential parts are again
among $0$, $z$, $Q_1$, $Q_2$, $Q_3$. Hence, for $r$ large and all $\theta$,
\[
|Y(re^{i\theta})|\le M\,r^K\exp\bigl(Q(r)-c\,r^{4/3}\bigr).
\]
Let $y_n=x_n-\lambda g_n=n![z^n]Y(z)$. Cauchy's inequality on the circle $|z|=z_n$ gives
$|y_n|\le n!\,Mz_n^K\exp\bigl(Q(z_n)-cz_n^{4/3}\bigr)z_n^{-n}$. On the other hand, Hayman's theorem
and \eqref{H-asym} give $g_n\ge c'\,n!\,z_n^{-11/6}e^{Q(z_n)}z_n^{-n}$ for $n$ large. Thus
$|y_n|/g_n=O\bigl(z_n^{K+2}e^{-cz_n^{4/3}}\bigr)\to0$, i.e. $x_n/g_n\to\lambda$.
\epf

\bl\label{lem:fupper}
Let $K_*=2^{3/10}$. There is a constant $M_*$ such that $f_n\le M_*K_*^n(n!)^{2/5}$ for all
$n\ge0$.
\el

\bpf
For $K>0$ let $w_n(K)=f_n/(K^n(n!)^{2/5})$. Dividing \eqref{posformF} by $K^n(n!)^{2/5}$ and using
$(n-s)!/n!=1/(n)_s$, we get
\begin{equation}\label{wineq}
w_n\le\frac{a_n}{K^n(n!)^{2/5}}+\sum_{s\ge1}\lambda_s(n,K)\,w_{n-s}
+\frac1K\sum_{i=2}^{n-3}\Bigl(n\binom{n-1}{i}\Bigr)^{-2/5}w_iw_{n-1-i}.
\end{equation}
Here $\lambda_1=K^{-1}n^{-2/5}$, $\lambda_2=0$, and for $s\ge3$, $\lambda_s(n,K)$ is the sum of
$\beta_{s+1}(n-s)K^{-s}(n)_s^{-2/5}$ and, if $s\ge5$, of
$\gamma_{s+2}(n-s)(n-s-1)K^{-s}(n)_s^{-2/5}$. We have $\beta_j\le j^2$ (indeed $\beta_4=\beta_5=2$
and $\beta_j=2+j(j-5)$ for $j\ge6$) and $\gamma_j\le2j$. Moreover, for $3\le s\le n$ we have
$(n)_s\ge n^3/9$ and, for $5\le s\le n$, $(n)_s\ge n^5/27$. Hence, uniformly in $n$,
\[
\lambda_s(n,K)\le\bar\lambda_s(K):=\bigl(3(s+1)^2+8(s+2)\bigr)K^{-s}\qquad(s\ge3).
\]
Furthermore, for fixed $s$ we have $\lambda_s(n,K)\to0$ as $n\to\infty$, except for $s=5$, where
the term $\gamma_7=2$ gives $\lambda_5(n,K)\to2K^{-5}$. For the quadratic term, note that
$\binom{n-1}i\ge\binom{n-1}2\ge n^2/6$ for $2\le i\le n-3$ and $n\ge5$, so the weights in
\eqref{wineq} sum to at most $3$. Moreover $\binom{n-1}{i}\ge\bigl(\frac{n-1}{i}\bigr)^i\ge2^i$
for $i\le(n-1)/2$.

\emph{Crude bound with $K=100$.} We show $w_n(100)\le4$ by induction. This holds for $n\le4$, since
$f_n\le4$. For $n\ge5$, \eqref{wineq} and the induction hypothesis give
\[
w_n\le n^2\,100^{-n}+4\Bigl(\frac1{100}+\sum_{s\ge3}\bar\lambda_s(100)\Bigr)+\frac{3\cdot16}{100}<1 .
\]

\emph{Descent.} Assume that $w_n(K)\le M$ for all $n$, and let
\[
2^{1/5}<K'\le K\le2^{1/5}K',\qquad\rho=K/K'\le2^{1/5}.
\]
Then $w_n(K')=\rho^nw_n(K)\le M\rho^n$. Let
$V_n=\max\bigl(1,\max_{k\le n}w_k(K')\bigr)$. In \eqref{wineq} for $K'$:
\begin{itemize}
\item the linear part is at most $V_{n-1}\sum_s\lambda_s(n,K')$, and this sum tends to
$2K'^{-5}<1$ by dominated convergence (with dominating sequence $\bar\lambda_s(K')$);
\item by symmetry, and bounding the factor with the smaller index by $M\rho^i$, the quadratic part
is at most $2K'^{-1}MV_{n-1}\,n^{-2/5}\sum_{i\ge2}(\rho2^{-2/5})^i$, which is $o(V_{n-1})$ since
$\rho2^{-2/5}\le2^{-1/5}$;
\item the first term tends to $0$ and $V_{n-1}\ge1$.
\end{itemize}
Hence $w_n(K')\le V_{n-1}$ for all $n\ge n_0$, so $(w_n(K'))$ is bounded.

Starting from $K=100$ and applying the descent with $K'=\max(2^{-1/5}K,K_*)$, we reach $K_*$ after
finitely many steps. Since $2^{1/5}<K_*$, every step is admissible.
\epf

\bp\label{olivier}
The limit $\kappa=\lim_{n\to\infty}f_n/g_n$ exists, and $1\le\kappa<\infty$.
\ep

\bpf
\emph{Boundedness of $f_n/g_n$.} Let $\Lambda_n=\max\bigl(1,\max_{2\le k\le n}f_k/g_k\bigr)$;
recall $f_0=f_1=0$. Since $T$ only involves indices $<n$ and has nonnegative coefficients,
\eqref{posformF} gives $f_n\le\Lambda_{n-1}(a+Tg)_n+q_n=\Lambda_{n-1}g_n+q_n$. By symmetry,
\[
q_n\le2\sum_{i=2}^{\lfloor(n-1)/2\rfloor}f_if_{n-1-i}\le\Lambda_{n-1}\,\kappa_n\,g_n,\qquad
\kappa_n:=\frac2{g_n}\sum_{i=2}^{\lfloor(n-1)/2\rfloor}f_i\,g_{n-1-i}.
\]
We show that $\kappa_n=O(n^{-6/5})$. By Stirling's formula,
\[
\Phi(n)=(2\pi)^{-1/5}(n!)^{2/5}2^{n/5}n^{-4/5}e^{\Pi(n)}(1+o(1)),
\]
where $\Pi(n)$ denotes the
argument of the exponential in the definition of $\Phi$; $\Pi$ is increasing for $n$ large.
Hence, by Corollary~\ref{cor:gasym}, there are $C$ and $n_1$ such that for $n\ge n_1$ and
$1\le k\le n/2+1$,
\[
\frac{g_{n-k}}{g_n}\le C\left(\frac{(n-k)!}{n!}\right)^{2/5}2^{-k/5}.
\]
Combined with Lemma~\ref{lem:fupper} (with $k=i+1$), this gives
\[
\frac{f_i\,g_{n-1-i}}{g_n}\le CM_*\Bigl(\frac{K_*}{2^{1/5}}\Bigr)^i\Bigl(n\binom{n-1}{i}\Bigr)^{-2/5}
\le CM_*\,n^{-2/5}\,t_i,
\]
where $t_i=\bigl(\rho'(i/(n-1))^{2/5}\bigr)^i$ and $\rho'=K_*2^{-1/5}=2^{1/10}$. Since
$t_{i+1}/t_i\le\rho'e^{2/5}\bigl(\frac{i+1}{n-1}\bigr)^{2/5}\le\frac12$
for $i+1\le\sqrt n$ and $n$ large, we have $\sum_{2\le i\le\sqrt n}t_i\le2t_2=O(n^{-4/5})$. For
$\sqrt n<i\le(n-1)/2$ we have $t_i\le(\rho'2^{-2/5})^i=2^{-3i/10}$. Altogether
$\kappa_n=O(n^{-2/5}\cdot n^{-4/5})=O(n^{-6/5})$.

Consequently $\Lambda_n\le\Lambda_{n-1}(1+\kappa_n)$, so $\Lambda_\infty=\sup_n\Lambda_n<\infty$,
and moreover $\sum_mq_m/g_m\le\sum_m\Lambda_\infty\kappa_m<\infty$.

\emph{Convergence.} By Lemma~\ref{lem:green}(i),
\[
\frac{f_n}{g_n}=1+\sum_{m=5}^nq_m\,\frac{e^{(m)}_n}{g_n}.
\]
By Lemma~\ref{lem:green}(iii) and Lemma~\ref{lem:limit}, $e^{(m)}_n/g_n\to\eps_m$ for each $m$, and
by Lemma~\ref{lem:green}(ii) the $m$-th term is bounded by $q_m/g_m$, which is summable. By
dominated convergence, $f_n/g_n\to\kappa=1+\sum_{m\ge5}q_m\eps_m\in[1,\infty)$.
\epf

\brem
The convergence of $f_n/g_n$ is slow (see Figure~\ref{fig1}). This comes from the tail
\[
\sum_{m>n}q_m\eps_m\le\sum_{m>n}\frac{q_m}{g_m}=O(n^{-1/5}).
\]
\erem

\begin{proof}[Proof of Theorem~\ref{thm:main2}]
The right-hand side of \eqref{asymBCK} is $A\Phi(n)$. By Corollary~\ref{cor:gasym},
$g_n\sim A_g\Phi(n)$, and by Proposition~\ref{olivier}, $f_n\sim\kappa g_n$. Hence
Theorem~\ref{thm:main2} holds with $A=\kappa A_g$.
\end{proof}

\begin{figure}[htb]
\begin{center}
\includegraphics[width=0.55\textwidth]{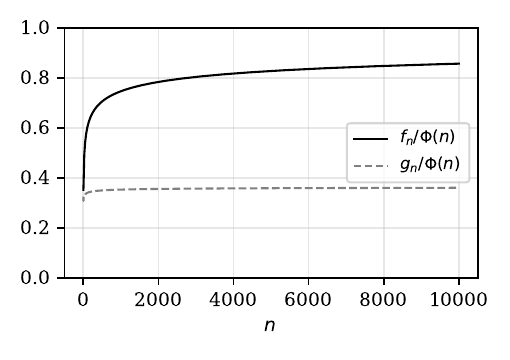}
\end{center}
\caption{\small The quotients $f_n/\Phi(n)$ (solid) and $g_n/\Phi(n)$ (dashed) for $n\le10000$,
where $\Phi(n)$ is the right-hand side of \eqref{asymBCK} without the constant $A$.}
\label{fig:2}
\end{figure}

We are unable to find the constant $A$ of Theorem~\ref{thm:main2} in closed form. A comparison
between $f_n$, $g_n$ and $\Phi(n)$ for the first 10000 terms of the sequences is given in
Figure~\ref{fig:2}.

\section{General $p$}\label{sec:general}

In this section $p\ge1$ is arbitrary, and $\mu_p,\nu_p,\varrho_p,\tau_p$ are as in the
introduction. Recall from Section~\ref{sec:prelim} the generating functions $F_p$ and $G_p$, the
coefficients $f^{(p)}_n$ and $g^{(p)}_n$, and Lemma~\ref{lem:linp}. Let
$H_p(z)=\sum_ng^{(p)}_nz^n/n!$ be the Borel transform of $G_p$. Writing $b_{k,j}=[z^j]\Gamma_k\ge0$,
we get as in \eqref{posform} the positive form
\begin{equation}\label{posformp}
g^{(p)}_n=a^{(p)}_n+\sum_{k=0}^p\sum_{j\ge3k+1}b_{k,j}\,(n-j+1)\cdots(n-j+k)\,g^{(p)}_{n-j+k},
\end{equation}
where $a^{(p)}_n=[z^n]\mathcal A_p$, and the same with $+\sum_{i}f^{(p)}_if^{(p)}_{n-1-i}$ for
$f^{(p)}$. Here $b_{p,3p+1}=2^p/p!$, and by Lemma~\ref{lem:linp}(ii) the coefficients $b_{k,j}$
grow at most polynomially in $j$.

\subsection{Growth, differential equation and the Borel transform}

\bp\label{prop:growthp}
Let $p\ge1$. There are constants $c,\kappa>0$ and, for every $\eps>0$, a constant $C_\eps$ such
that for $n\ge2$
\[
c\,n^{-\kappa}(n!)^{\mu_p}\varrho_p^{\,n}\le g^{(p)}_n\le f^{(p)}_n\le
C_\eps\,(n!)^{\mu_p}\bigl((1+\eps)\varrho_p\bigr)^n .
\]
In particular $\log f^{(p)}_n=\mu_p\log n!+n\log\varrho_p+o(n)$, the same holds for
$g^{(p)}_n$, and $\log H_p(r)\ge\tau_pr^{\nu_p}-O(\log r)$ as $r\to+\infty$.
\ep

\bpf
Write $\mu=\mu_p$ and $\varrho=\varrho_p$.

\emph{Lower bound.} Keeping only the term $k=p$, $j=3p+1$ in \eqref{posformp} gives, for
$n\ge3p+1$,
\[
g_n\ge\frac{2^p}{p!}(n-3p)\cdots(n-2p-1)\,g_{n-2p-1}\ge\frac{2^p}{p!}(n-3p)^p\,g_{n-2p-1}.
\]
Let $m_j=n-j(2p+1)$. Let $J$ be the first index with $m_J<6p$; thus $4p-1\le m_J<6p$ and
$m_j\ge6p$ for $j<J$. Iterating,
$g_n\ge(2^p/p!)^J\prod_{j<J}(m_j-3p)^p\,g_{m_J}$. On the other hand, since
$\varrho^{2p+1}=2^p/p!$ and $(m)_{2p+1}^{\mu}\le m^{(2p+1)\mu}=m^p$,
\[
(n!)^\mu\varrho^n=(m_J!)^\mu\varrho^{m_J}\prod_{j<J}\bigl((m_j)_{2p+1}\bigr)^\mu\varrho^{2p+1}
\le(m_J!)^\mu\varrho^{m_J}\Bigl(\frac{2^p}{p!}\Bigr)^J\prod_{j<J}m_j^p .
\]
Let $c_0=\min\{g_m/((m!)^\mu\varrho^m):4p-1\le m<6p\}$; then $c_0>0$, since $g_m\ge1$ for $m\ge2$.
Hence
\[
\frac{g_n}{(n!)^\mu\varrho^n}\ge c_0\prod_{j<J}\Bigl(1-\frac{3p}{m_j}\Bigr)^p
\ge c_0\exp\Bigl(-6p^2\sum_{j<J}\frac1{m_j}\Bigr),
\]
where we used $\log(1-x)\ge-2x$ for $0\le x=3p/m_j\le\frac12$. Since the $m_j$, $j<J$, form an
arithmetic progression with difference $2p+1$ contained in $[6p,n]$,
$\sum_{j<J}1/m_j\le\frac1{6p}+\frac1{2p+1}\log n$. This gives the lower bound with
$\kappa=6p^2/(2p+1)$. Finally $f\ge g$, since $zF_p^2$ has nonnegative coefficients.

\emph{Upper bound: the normalised inequality.} For $K>0$ let $w_n(K)=f_n/(K^n(n!)^\mu)$. Dividing
the positive form of $f$ by $K^n(n!)^\mu$ gives
\begin{equation}\label{wineqp}
w_n\le\frac{a^{(p)}_n}{K^n(n!)^\mu}+\sum_{k,j}\lambda_{k,j}(n,K)\,w_{n-j+k}
+\frac1K\sum_{i=2}^{n-3}\Bigl(n\binom{n-1}i\Bigr)^{-\mu}w_iw_{n-1-i},
\end{equation}
where, with $s=j-k\ge2k+1$,
$\lambda_{k,j}(n,K)=b_{k,j}\,(n-j+1)\cdots(n-j+k)\,K^{-s}(n)_s^{-\mu}$ for $j\le n$, and
$\lambda_{k,j}(n,K)=0$ otherwise. We use four estimates.
\begin{itemize}
\item[(a)] \emph{Uniform majorant.} We have $b_{k,j}\le Bj^E$ and
$(n-j+1)\cdots(n-j+k)\le n^k$. Moreover $(n)_s=(n)_{2k+1}(n-2k-1)_{s-2k-1}\ge(n)_{2k+1}(s-2k-1)!$.
For $n\ge2k+1$, $(n)_{2k+1}\ge(n/(2k+1))^{2k+1}$ and $(2k+1)\mu\ge k$, hence
$n^k(n)_{2k+1}^{-\mu}\le C_0$. Therefore
$\lambda_{k,j}(n,K)\le\bar\lambda_{k,j}(K):=C_0Bj^EK^{-s}((s-2k-1)!)^{-\mu}$ for all $n$. For every
$K>0$ we have $\sum_{k,j}\bar\lambda_{k,j}(K)<\infty$, also when $K<1$ (which happens for $p\ge4$,
since then $\varrho_p<1$), and this sum tends to $0$ as $K\to\infty$.
\item[(b)] \emph{Limits.} For fixed $(k,j)$, $\lambda_{k,j}(n,K)$ is of order $n^{k-\mu s}$, and
$k-\mu s\le k-\mu(2k+1)=\frac{k-p}{2p+1}\le0$, with equality only for $k=p$, $s=2p+1$. Hence
$\lambda_{k,j}(n,K)\to0$ except for $(k,j)=(p,3p+1)$, where the limit is
$\frac{2^p}{p!}K^{-(2p+1)}$. By (a) and dominated convergence,
$\sum_{k,j}\lambda_{k,j}(n,K)\to\frac{2^p}{p!}K^{-(2p+1)}$.
\item[(c)] \emph{Quadratic weights.} For $n\ge5$ and $2\le i\le n-3$,
$\binom{n-1}i\ge\binom{n-1}2\ge n^2/6$. Hence the weights in \eqref{wineqp} sum to at most
$n(n^3/6)^{-\mu}\le6$, because $3\mu\ge1$. For $i\le(n-1)/2$ we have $\binom{n-1}i\ge2^i$, so the
weight is at most $n^{-\mu}2^{-\mu i}$.
\item[(d)] \emph{Source term.} $a^{(p)}_n$ grows polynomially, since $\mathcal A_p$ is rational with
its only pole at $z=1$.
\end{itemize}

\emph{Upper bound: crude bound.} Let $W=\max(3,f_0,\dots,f_4)$, and choose $K\ge\max(1,18W)$
such that $a^{(p)}_n\le K^n(n!)^\mu$ for all $n$ and $\sum_{k,j}\bar\lambda_{k,j}(K)\le\frac13$. Then
$w_n(K)\le f_n\le W$ for $n\le4$. For $n\ge5$, assuming $w_m(K)\le W$ for $m<n$, \eqref{wineqp}
gives
\[
w_n(K)\le1+\frac W3+\frac{6W^2}K\le1+\frac{2W}3\le W.
\]

\emph{Upper bound: descent.} Assume $w_n(K)\le M$ for all $n$, and let
$\varrho<K'\le K\le2^{\mu/2}K'$. Put $\rho=K/K'\le2^{\mu/2}$; then
$w_n(K')=\rho^nw_n(K)\le M\rho^n$. Let $V_n=\max\bigl(1,\max_{m\le n}w_m(K')\bigr)$, and write
$\frac{2^p}{p!}K'^{-(2p+1)}=1-2\eta$ with $\eta>0$. In \eqref{wineqp} for $K'$, we bound the
factor with the smaller index in each quadratic term by $M\rho^i$, and use (c). This gives
\[
w_n(K')\le\frac{a^{(p)}_n}{K'^n(n!)^\mu}+V_{n-1}\sum_{k,j}\lambda_{k,j}(n,K')
+\frac{2M}{K'}V_{n-1}\,n^{-\mu}\sum_{i\ge2}(\rho2^{-\mu})^i .
\]
By (b), the middle sum tends to $1-2\eta$. Since $\rho2^{-\mu}\le2^{-\mu/2}<1$, the last term is at
most $\frac\eta2V_{n-1}$ for $n$ large. By (d) and $V_{n-1}\ge1$, so is the first term. Hence
$w_n(K')\le V_{n-1}$ for $n\ge n_0$, and $(w_n(K'))$ is bounded.

Given $\eps>0$, we start from the crude bound and apply the descent with
$K'=\max\bigl(2^{-\mu/2}K,(1+\eps)\varrho\bigr)$. After finitely many steps we reach
$K'=(1+\eps)\varrho$, which proves the upper bound.

\emph{Borel transform.} Take $n=\lfloor(\varrho_pr)^{\nu_p}\rfloor$. By the lower bound and
Stirling's formula,
$\log\bigl(g_nr^n/n!\bigr)\ge n\log(\varrho_pr)-(1-\mu_p)(n\log n-n)-O(\log n)
=(1-\mu_p)(\varrho_pr)^{\nu_p}-O(\log r)$, since $(1-\mu_p)\nu_p=1$.
\epf

\bp\label{prop:odep}
Let $q_k(x)=(1-x)^{2p-1}\Gamma_k(x)/x^k$ for $0\le k\le p$. Then $q_k$ is a polynomial of degree
at most $3p$, divisible by $x^{2k+1}$, and its coefficient of $x^{2p+1}$ for $k=p$ is $2^p/p!$.
The function $H_p$ satisfies $L_pH_p=0$, where $D=\mathrm d/\mathrm dz$ and
\[
L_p=D^{p+1}(D-1)^{2p-1}-\sum_{k=0}^p\Pi_k(D)\,z^kD^k,\qquad
\Pi_k(D)=\sum_mq_{k,m}D^{3p-m}.
\]
Thus $H_p$ satisfies a linear differential equation of order $3p$ whose coefficients are
polynomials of degree at most $p$.
\ep

\bpf
The first assertions follow from Lemma~\ref{lem:linp}. Let $\mathfrak B$ denote the Borel
transform and $J$ the integration operator $(JX)(z)=\int_0^zX$. Then $\mathfrak B(zX)=J\mathfrak BX$,
and $\mathfrak B$ commutes with $\theta=zD$. Since $z^kD^k=\theta(\theta-1)\cdots(\theta-k+1)$, we
have $\mathfrak B(z^kG^{(k)})=z^kD^k\mathfrak BG$, hence $\mathfrak B(\Gamma_kG^{(k)})=
c_k(J)\,z^kD^kH_p$ with $c_k(x)=\Gamma_k(x)/x^k$. Applying $\mathfrak B$ to the linear equation
of Lemma~\ref{lem:linp} and multiplying by $(1-J)^{2p-1}$ yields
\[
(1-J)^{2p-1}H_p=\mathfrak B\bigl((1-z)^{2p-1}\mathcal A_p\bigr)+\sum_kq_k(J)\,z^kD^kH_p .
\]
The first term on the right is a polynomial of degree at most $2p$. We now apply $D^{3p}$. Since
$DJ=I$, we have $D^{3p}J^m=D^{3p-m}$ for $m\le3p$ and $D^{2p-1}(1-J)^{2p-1}=(D-1)^{2p-1}$, and
$D^{3p}$ annihilates the polynomial.
\epf

\brem
For $p=2$, Lemma~\ref{lem:diffeq} gives an equation of order $5$ for $H$, smaller than the bound
$3p=6$. For $p=1$, the operator $L_1=D^3-2D^2+D-2zD$ is already of order
$3=2p+1$. Writing $H_1'=e^zu$, it reduces to the Airy equation $u''=2zu$, so
$H_1'(z)=e^z\bigl(a\operatorname{Ai}(2^{1/3}z)+b\operatorname{Bi}(2^{1/3}z)\bigr)$. Whether the
minimal order is $2p+1$ for every $p$ remains open; Proposition~\ref{prop:odep} only gives the
bound $3p$.
\erem

\begin{theorem}\label{thm:Hp}
Let $p\ge1$. There exist $\delta>0$, $C>0$, a real number $\sigma_p$ and real numbers
$a_1,\dots,a_{2p}$ such that, uniformly as $z\to\infty$ in the sector $S$ of the form
\eqref{sector},
\[
H_p(z)=C\,z^{\sigma_p}\exp\Bigl(\tau_p\,z^{(2p+1)/(p+1)}+\sum_{m=1}^{2p}a_mz^{m/(p+1)}\Bigr)
\bigl(1+O(|z|^{-1/(p+1)})\bigr).
\]
\end{theorem}

\bpf
\emph{The Newton polygon.} Write $L_p=\sum_iP_i(z)D^i$. Then $P_{3p}=1$ and
$\deg P_i\le p$. The monomials coming from $\Pi_k(D)z^kD^k$ are obtained by normal ordering
$D^az^k=\sum_j\binom aj\frac{k!}{(k-j)!}z^{k-j}D^{a-j}$ with $a=3p-m$ and $m\ge2k+1$. They are
of the form $z^dD^i$ with $d=k-j$ and $3p-i=m-k+j\ge k+1+j$. Hence
\[
d\le\frac p{p+1}(3p-i),
\]
with equality only for $j=0$, $k=p$, $m=2p+1$, i.e. for the monomial
$-\frac{2^p}{p!}z^pD^{2p-1}$. Indeed $(k+1)\frac p{p+1}\ge k$, with equality iff $k=p$. So the
upper edge of the Newton polygon joins $(3p,0)$ to $(2p-1,p)$, has slope $p/(p+1)$, contains no
other monomial, and its characteristic equation is $\lambda^{p+1}=2^p/p!$. All other edges
have smaller slopes.

\emph{Formal solutions.} For each of the $p+1$ (simple) roots $\lambda$ we construct a formal
solution as in the proof of Proposition~\ref{fracpower}. Writing
$y=\exp\int w$, we have $L_py=y\,\mathcal R(w)$ for a differential polynomial $\mathcal R$, and we
solve $\mathcal R(w)=0$ with $w=\sum_{j\ge0}w_jz^{(p-j)/(p+1)}$ and $w_0=\lambda$.
\begin{itemize}
\item The top coefficient of $\mathcal R(w)$ is
$w_0^{2p-1}\bigl(w_0^{p+1}-2^p/p!\bigr)=0$.
\item By the Newton polygon analysis above, the coefficient determining $w_j$ is
\[
\Bigl(3p\,w_0^{3p-1}-(2p-1)\frac{2^p}{p!}w_0^{2p-2}\Bigr)w_j+(\text{terms in }w_0,\dots,w_{j-1}),
\]
and the factor in front of $w_j$ equals $(p+1)\frac{2^p}{p!}w_0^{2p-2}\ne0$.
\item Integrating, the term $j=2p+1$ produces the power $z^{\sigma_\lambda}$ with
$\sigma_\lambda=w_{2p+1}$, and the terms $j\le2p$ produce
$Q_\lambda(z)=\sum_{j=0}^{2p}\frac{p+1}{2p+1-j}w_jz^{(2p+1-j)/(p+1)}$.
\end{itemize}
For $\lambda_0=(2^p/p!)^{1/(p+1)}>0$ all $w_j$ are real, and the leading term of
$Q_{\lambda_0}$ is $\tau_pz^{\nu_p}$. The upper edge of the Newton polygon has slope
$\gamma_{\max}=p/(p+1)$ and length $p+1$. By (F1) of Section~\ref{sec:HT}, exactly $p+1$ elements
of a formal fundamental system have exponential parts of degree $1+\gamma_{\max}=\nu_p$. The $p+1$
formal solutions just constructed have distinct exponential parts of this degree, so they are
linearly independent and can be taken as these elements. All other formal solutions have
exponential parts that vanish or have some degree $\nu'<\nu_p$; their real parts are
$O(|z|^{\nu'})$.

\emph{Actual solutions and domination.} By (F2), in a small sector $S$ there is a fundamental system $(Y_\lambda)\cup(Y'_i)$ asymptotic to these formal
solutions, and $H_p=\sum_\lambda c_\lambda Y_\lambda+\sum_ic'_iY'_i$. For $\lambda\ne\lambda_0$ we
have $\arg\lambda=2\pi k/(p+1)$ with $k\ne0$, so for $|\arg z|<\delta$ small
\[
\operatorname{Re}Q_\lambda(z)\le\cos\bigl(\tfrac{2\pi}{p+1}-\nu_p\delta\bigr)\,\tau_p|z|^{\nu_p}+O(|z|^{2p/(p+1)}),
\]
while $\operatorname{Re}Q_{\lambda_0}(z)\ge\cos(\nu_p\delta)\,\tau_p|z|^{\nu_p}-O(|z|^{2p/(p+1)})$.
The other solutions grow at most like $\exp(O(|z|^{\nu'}))$. Hence, if $c_{\lambda_0}\ne0$,
$H_p=c_{\lambda_0}Y_{\lambda_0}\bigl(1+O(e^{-\eta|z|^{\nu_p}})\bigr)$ in $S$, which is the
assertion.

\emph{Non-vanishing of $c_{\lambda_0}$ and positivity of $C$.} Suppose $c_{\lambda_0}=0$. On the
positive axis, $H_p$ is then a combination of the $Y_\lambda$ with $\lambda\ne\lambda_0$ and of the
$Y'_i$. For $\lambda=\lambda_0e^{2\pi ik/(p+1)}$ with $1\le k\le p$, we have
$\operatorname{Re}Q_\lambda(r)=\cos\bigl(\frac{2\pi k}{p+1}\bigr)\tau_pr^{\nu_p}+O(r^{2p/(p+1)})$ and
$\cos\frac{2\pi k}{p+1}\le\cos\frac{2\pi}{p+1}$. The $Y'_i$ are $O(\exp(O(r^{\nu'})))$ with
$\nu'<\nu_p$. They may grow (for instance like $e^r$), which is why the bound below involves the
maximum with $0$. Hence
\[
\log|H_p(r)|\le\max\Bigl\{0,\cos\Bigl(\frac{2\pi}{p+1}\Bigr)\Bigr\}\tau_pr^{\nu_p}+o(r^{\nu_p}).
\]
Since $\max\{0,\cos\frac{2\pi}{p+1}\}<1$, this contradicts the lower bound
$\log H_p(r)\ge\tau_pr^{\nu_p}-O(\log r)$ of Proposition~\ref{prop:growthp}. Finally,
\[
C=\lim_{r\to\infty}H_p(r)\,r^{-\sigma_p}e^{-Q_{\lambda_0}(r)}
\]
is a limit of positive numbers and is nonzero, so $C>0$.
\epf

\brem\label{rem:conjp}
\begin{enumerate}
\item[(a)] The leading constant is $\tau_p=\frac{p+1}{2p+1}(2^p/p!)^{1/(p+1)}$. The value
$\frac{p+1}{2^p+1}(2^p/p!)^{1/(p+1)}$ conjectured in \cite{BG14} coincides with
$\tau_p$ only for $p=1,2$.
\item[(b)] The coefficients $a_m$ and $\sigma_p$ are obtained from the triangular recursion in the
proof of Theorem~\ref{thm:Hp}, which determines one coefficient per linear equation.
\item[(c)] Theorem~\ref{thm:Hp} is used in Section~\ref{sec:fp} below to obtain the
asymptotic number of $BCK(p)$-terms for every $p\ge2$.
\end{enumerate}
\erem

\subsection{The asymptotic number of $BCK(p)$-terms for $p\ge2$}\label{sec:fp}

In this subsection, $Q=Q_{\lambda_0}$, $\sigma_p$, $C$ and $a_1,\dots,a_{2p}$ are as in
Theorem~\ref{thm:Hp}, and we write $a_{2p+1}=\tau_p$, so that
$Q(z)=\sum_{m=1}^{2p+1}a_mz^{m/(p+1)}$.

\bl\label{lem:a2p}
For $p\ge2$ we have $a_{2p}=\frac14\bigl(2^p/p!\bigr)^{2/(p+1)}>0$. For $p=1$, $a_2=1$.
\el

\bpf
Give the monomial $z^dD^i$ the weight $\varpi=(p+1)d+pi$. By the Newton polygon analysis in the
proof of Theorem~\ref{thm:Hp}, the two monomials of the upper edge have weight $3p^2$ and all other
monomials of $L_p$ have weight at most $3p^2-1$.
\begin{itemize}
\item A monomial of weight $3p^2-1$ must have $d\equiv-1\pmod p$, i.e. $d=p-1$ and $i=2p$ (we use
$p\ge2$).
\item In the normal ordering of $\Pi_k(D)z^kD^k$, the monomial $z^{p-1}D^{2p}$ arises only from
$k=p-1$, $j=0$ and $m=2p-1=2k+1$. Its coefficient in $L_p$ is therefore
$-q_{p-1,2p-1}=-[z^{3p-2}]\Gamma_{p-1}=-2^{p-1}/(p-1)!$. Indeed, by the proof of
Lemma~\ref{lem:linp}, the lowest coefficient of $\Gamma_k$ is $2^k/k!$, coming from $l=I=k$.
\end{itemize}
In the Riccati equation $\mathcal R(w)=0$ of that proof, a monomial $z^dD^i$ contributes $z^dw^i$ at
the exponent $\varpi/(p+1)$. Every term containing a derivative of $w$ has an exponent lower by
at least $(2p+1)/(p+1)$. Hence the coefficient of $z^{(3p^2-1)/(p+1)}$ in $\mathcal R(w)$ is
\[
(p+1)\frac{2^p}{p!}\,w_0^{2p-2}\,w_1-\frac{2^{p-1}}{(p-1)!}\,w_0^{2p}=0,
\qquad\text{so}\qquad w_1=\frac{p\,w_0^2}{2(p+1)}.
\]
Since $a_{2p}=\frac{p+1}{2p}w_1$ and $w_0^{p+1}=2^p/p!$, the claim follows. For $p=1$, the
representation by Airy functions gives the exponent $\tau_1z^{3/2}+z$.
\epf

For $p=2$ this gives $a_4=2^{-4/3}$, the coefficient $\alpha_4$ used in the proof of
Lemma~\ref{lem:hayman}.

\bl\label{lem:haymanp}
For every $p\ge1$, the function $H_p$ is Hayman-admissible, with $\eta(r)=r^{-5\nu_p/12}$ and
$b(r)\sim\nu_p^2\tau_pr^{\nu_p}$.
\el

\bpf
We follow the proof of Lemma~\ref{lem:hayman}; parts (a)--(c) carry over verbatim, with the
following changes.
\begin{itemize}
\item \emph{Taylor estimate.} The third derivative of $\phi\mapsto\log H_p(re^{i\phi})$ is
$O(r^{\nu_p})$, and $r^{\nu_p}\eta(r)^3\to0$ while $b(r)\eta(r)^2\to\infty$.
\item \emph{Arc inside the sector.} The terms with $a_m<0$ (necessarily $m\le2p-1$) contribute at
most $O(r^{(2p-1)/(p+1)}\phi^2)$. Hence $\operatorname{Re}Q(re^{i\phi})-Q(r)\le-c\,r^{\nu_p}\phi^2$
for $|\phi|\le\delta$.
\end{itemize}
For part~(d), which rests on (F3), note that $L_p$ has polynomial coefficients in $z=t^{p+1}$. The substitution
$t\mapsto e^{2\pi i/(p+1)}t$ therefore maps formal solutions to formal solutions and permutes the
roots of $\lambda^{p+1}=2^p/p!$. Since these roots are simple, the exponential parts of the $p+1$
dominant formal solutions are obtained from $Q$ by replacing $z^{1/(p+1)}$ with its $p+1$
determinations. For $z=re^{i\theta}$, their real parts form the set $\{F_r(\psi):\psi\in\Psi(\theta)\}$,
where
\[
F_r(\psi)=\sum_{m=1}^{2p+1}a_mr^{m/(p+1)}\cos(m\psi),\qquad
\Psi(\theta)=\Bigl\{\tfrac{\theta+2\pi k}{p+1}:0\le k\le p\Bigr\}.
\]
The other formal solutions have exponential parts $O(r^{\nu'})$ with $\nu'<\nu_p$.

For $\delta\le|\theta|\le\pi$, every $\psi\in\Psi(\theta)$ satisfies, modulo $2\pi$,
$\delta/(p+1)\le|\psi|\le\pi$. The function $\max\bigl(1-\cos((2p+1)\psi),\,1-\cos(2p\psi)\bigr)$
vanishes only on $2\pi\mathbb Z$, so it is bounded below by some $\mu_\delta>0$ there. Using
$\tau_p>0$ and $a_{2p}>0$ (Lemma~\ref{lem:a2p}), we get
\[
F_r(\psi)\le Q(r)-\min(\tau_p,a_{2p})\,\mu_\delta\,r^{2p/(p+1)}+O(r^{(2p-1)/(p+1)}).
\]
The rest of part~(d) is unchanged, and yields
$|H_p(re^{i\theta})|/H_p(r)=O\bigl(r^{K}e^{-c\,r^{2p/(p+1)}}\bigr)$ uniformly for $\theta\in A$.
\epf

Let $\tilde\Lambda(r)=\log C+\sigma_p\log r+Q(r)$. For real $x>0$, let $\tilde\zeta(x)$ be the
solution of $r\tilde\Lambda'(r)=x$, let $\tilde b(r)=r(r\tilde\Lambda'(r))'$, and put
\[
\Psi_p(x)=\Gamma(x+1)\,\frac{\exp\bigl(\tilde\Lambda(\tilde\zeta)-x\log\tilde\zeta\bigr)}
{\sqrt{2\pi\tilde b(\tilde\zeta)}},\qquad\tilde\zeta=\tilde\zeta(x).
\]

\bl\label{lem:ratiop}
For every $p\ge1$ we have $g^{(p)}_n\sim\Psi_p(n)$. Moreover there are $C'$ and $n_1$ such that
\[
\frac{g^{(p)}_{n-k}}{g^{(p)}_n}\le C'\left(\frac{(n-k)!}{n!}\right)^{\mu_p}\varrho_p^{-k}
\qquad\text{whenever }n-k\ge n_1 .
\]
\el

\bpf
By Hayman's theorem and Lemma~\ref{lem:haymanp},
$g_n\sim n!\,H_p(\zeta_n)\zeta_n^{-n}/\sqrt{2\pi b(\zeta_n)}$. On the positive axis we have
$\log H_p=\tilde\Lambda+O(r^{-1/(p+1)})$, and by Cauchy's estimates the same holds for the first
two logarithmic derivatives. Hence the minima of $\log H_p(r)-n\log r$ and of
$\tilde\Lambda(r)-n\log r$ differ by $o(1)$, $\zeta_n/\tilde\zeta(n)\to1$ and
$b(\zeta_n)/\tilde b(\tilde\zeta(n))\to1$. This proves $g_n\sim\Psi_p(n)$.

Now let $\ell(x)=\log\Psi_p(x)-\mu_p\log\Gamma(x+1)-x\log\varrho_p$. By the envelope theorem,
\[
\ell'(x)=(1-\mu_p)\,\psi(x+1)-\log\tilde\zeta(x)-\log\varrho_p+O(1/x),
\]
where $\psi$ is the digamma function. Here $1-\mu_p=1/\nu_p$ and $\psi(x+1)=\log x+O(1/x)$.
Moreover
\[
x=r\tilde\Lambda'(r)=\nu_p\tau_pr^{\nu_p}\bigl(1+\varkappa r^{-1/(p+1)}+O(r^{-2/(p+1)})\bigr),
\qquad r=\tilde\zeta(x),
\]
where $\varkappa=\frac{2p\,a_{2p}}{(p+1)\nu_p\tau_p}>0$.
Since $\nu_p\tau_p=\varrho_p^{\nu_p}$, this gives
\[
\ell'(x)=\frac{\varkappa}{\nu_p}\,\tilde\zeta(x)^{-1/(p+1)}+O\bigl(\tilde\zeta(x)^{-2/(p+1)}\bigr)+O(1/x).
\]
As $\tilde\zeta(x)^{-1/(p+1)}\asymp x^{-1/(2p+1)}$, we get $\ell'(x)>0$ for $x\ge x_0$. Hence, for
$n-k\ge x_0$,
\[
\log\frac{\Psi_p(n-k)}{\Psi_p(n)}=\mu_p\log\frac{(n-k)!}{n!}-k\log\varrho_p+\ell(n-k)-\ell(n)
\le\mu_p\log\frac{(n-k)!}{n!}-k\log\varrho_p,
\]
and $g_n\sim\Psi_p(n)$ gives the claim.
\epf

For a sequence $x$, let $R^{(p)}_n(x)$ denote the coefficient of $z^n$ in
$(1-z)^{2p-1}\bigl(X(z)-\sum_k\Gamma_k(z)X^{(k)}(z)\bigr)$, where $X(z)=\sum_nx_nz^n$. By
Lemma~\ref{lem:linp}(ii), this is a linear recurrence with polynomial coefficients, of finite
order. The sequence $g^{(p)}$ satisfies $R^{(p)}_n(g^{(p)})=0$ for $n>2p$.

\bl\label{lem:limitp}
Let $p\ge1$. If $R^{(p)}_n(x)=0$ for all $n\ge N$, then $x_n/g^{(p)}_n$ converges to a finite
limit.
\el

\bpf
This is the proof of Lemma~\ref{lem:limit}, with $L_7$ replaced by $D^s\circ L_p$ for $s$ large.
\begin{itemize}
\item The recurrence $R^{(p)}_n(x)=0$ reads $x_n=\sum_{s=1}^S\pi_s(n)x_{n-s}$ with
$\deg\pi_s\le\mu_ps$.
\begin{itemize}
\item The factor $(1-z)^{2p-1}$ alone contributes constant coefficients.
\item The term $\Gamma_kX^{(k)}$, multiplied by the polynomial $(1-z)^{2p-1}\Gamma_k$ of valuation
$3k+1$, contributes coefficients of degree $k$ at shifts $s\ge2k+1$. If $s\le2p+1$, then
$k\le\frac{s-1}2\le\mu_ps$; if $s>2p+1$, then $k\le p<\mu_ps$.
\end{itemize}
By Lemma~\ref{lem:entire}, $X(z)=\sum_nx_nz^n/n!$ is entire.
\item As in the proof of Proposition~\ref{prop:odep}, the Borel transform $X$ satisfies
$(1-J)^{2p-1}(X-\sum_kc_k(J)z^kD^kX)=P$ for a polynomial $P$. Applying $D^{3p+s}$ with
$s>\deg P$ gives $D^sL_pX=0$.
\item Normal ordering $D^sP_i(z)$ only adds monomials of lower weight. So the upper edge of the
Newton polygon of $D^sL_p$ is the same as for $L_p$. By (F1), the additional formal solutions have
exponential parts of degree less than $\nu_p$.
\item In $S$, let $\lambda$ be such that $Y=X-\lambda H_p$ has no component on the dominant
solution. The domination argument in the proof of Theorem~\ref{thm:Hp} gives
$|Y|\le\exp\bigl(Q(r)-c\,r^{\nu_p}\bigr)$ in $S$.
\item Part~(d) of the proof of Lemma~\ref{lem:haymanp} relies on (F3). It therefore applies to the
solutions $X$ and $H_p$ of $D^sL_py=0$ and gives
\[
|Y(re^{i\theta})|\le Mr^K\exp\bigl(Q(r)-c\,r^{2p/(p+1)}\bigr)\quad\text{outside }S.
\]
\end{itemize}
By Cauchy's inequality at $r=\zeta_n$ and the lower bound
$g_n\ge c'\,n!\,\zeta_n^{\sigma_p-\nu_p/2}e^{Q(\zeta_n)}\zeta_n^{-n}$ (Hayman), the quotient
$|x_n-\lambda g_n|/g_n$ tends to $0$.
\epf

\bp\label{prop:olivierp}
Let $p\ge2$. Then $\kappa_p=\lim_{n\to\infty}f^{(p)}_n/g^{(p)}_n$ exists, and
$1\le\kappa_p<\infty$.
\ep

\bpf
The operator $(T_px)_n=\sum_{k,j}b_{k,j}(n-j+1)\cdots(n-j+k)\,x_{n-j+k}$ has nonnegative
coefficients and only involves indices $<n$. Moreover $g=a+T_pg$ and $f=a+T_pf+q$, where
$q_n=\sum_if_if_{n-1-i}$. Lemma~\ref{lem:green} holds verbatim for $T_p$, with
item~(iii) replaced by: $e^{(m)}$ satisfies $R^{(p)}_n(e^{(m)})=0$ for $n\ge m+2p$. Its proof
also carries over: the source $z^m(1-z)^{2p-1}$ has degree $m+2p-1$.

\emph{Boundedness.} As in the proof of Proposition~\ref{olivier},
$\Lambda_n\le\Lambda_{n-1}(1+\kappa_n)$ with
$\kappa_n=\frac2{g_n}\sum_{i=2}^{\lfloor(n-1)/2\rfloor}f_ig_{n-1-i}$. By
Proposition~\ref{prop:growthp} with $1+\eps=2^{\mu_p/2}$, we have
$f_i\le M\bigl(2^{\mu_p/2}\varrho_p\bigr)^i(i!)^{\mu_p}$. Combined with Lemma~\ref{lem:ratiop},
this gives
\[
\frac{f_ig_{n-1-i}}{g_n}\le C''\,n^{-\mu_p}\,t_i,\qquad
t_i=\Bigl(2^{\mu_p/2}\Bigl(\frac i{n-1}\Bigr)^{\mu_p}\Bigr)^i ,
\]
using $\binom{n-1}{i}\ge\bigl(\frac{n-1}i\bigr)^i$.
\begin{itemize}
\item For $i+1\le\sqrt n$ and $n$ large, $t_{i+1}/t_i\le\frac12$, so these terms contribute
$O(t_2)=O(n^{-2\mu_p})$.
\item For $i>\sqrt n$ we have $t_i\le2^{-\mu_pi/2}$, which is negligible.
\end{itemize}
Hence $\kappa_n=O(n^{-3\mu_p})$. This is summable because $3\mu_p=\frac{3p}{2p+1}>1$ for $p\ge2$.
Thus $\Lambda_n$ is bounded and $\sum_mq_m/g_m<\infty$.

\emph{Convergence.} By Lemma~\ref{lem:limitp}, each $e^{(m)}_n/g_n$ converges. Dominated
convergence then concludes exactly as in the proof of Proposition~\ref{olivier}.
\epf

\bp\label{thm:fp}
Let $p\ge2$. The number of $BCK(p)$-terms of size $n$ satisfies
\[
f^{(p)}_n\sim\kappa_p\,\Psi_p(n)
=\kappa_p\,n!\,\frac{C\,\tilde\zeta^{\,\sigma_p-n}\exp\bigl(Q(\tilde\zeta)\bigr)}{\sqrt{2\pi\tilde b(\tilde\zeta)}},
\qquad\tilde\zeta=\tilde\zeta(n).
\]
Expanding the saddle point gives the form
\[
f^{(p)}_n\sim A_p\,n^{\mu_pn}\bigl(\varrho_pe^{-\mu_p}\bigr)^n
\exp\Bigl(\sum_{m=1}^{2p}e_mn^{m/(2p+1)}\Bigr)\,n^{\sigma_p/\nu_p},
\]
with constants $e_m$ computable from $a_1,\dots,a_{2p+1}$ and $\sigma_p$. For $p=2$ this is
\eqref{asymBCK}.
\ep

\bpf
The first formula follows from Proposition~\ref{prop:olivierp} and Lemma~\ref{lem:ratiop}. The
second follows by inserting the expansion of $\tilde\zeta(n)$ in powers of $n^{-1/(2p+1)}$, as in
Corollary~\ref{cor:gasym}. The leading factors come from $(1-\mu_p)\nu_p=1$ and
$\nu_p\tau_p=\varrho_p^{\nu_p}$.
\epf

\begin{figure}[htb]
\begin{center}
\includegraphics[width=0.55\textwidth]{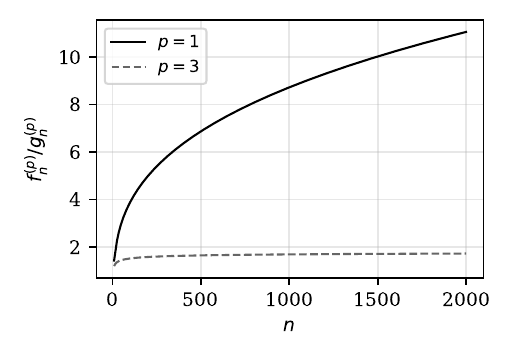}
\end{center}
\caption{\small The quotients $f^{(p)}_n/g^{(p)}_n$ for $p=1$ (solid) and $p=3$ (dashed),
$n\le2000$.}
\label{fig:p13}
\end{figure}

\brem\label{rem:p1}
For $p=1$ the proof of Proposition~\ref{prop:olivierp} breaks down: $3\mu_1=1$, and the bound
$\kappa_n=O(n^{-1})$ is no longer summable. This is not a technical limitation.
\begin{itemize}
\item Theorem~\ref{thm:Hp} and Lemmas~\ref{lem:haymanp} and~\ref{lem:ratiop} hold for $p=1$.
\item The Airy representation of $H_1'$ gives $Q(z)=\tau_1z^{3/2}+z$ with
$\tau_1=\frac{2\sqrt2}3$, and $\sigma_1=-\frac34$.
\item Expanding $\Psi_1(n)$ as in the proof of Corollary~\ref{cor:gasym} gives
\[
g^{(1)}_n\sim C'\,n^{-1/2}\Bigl(\frac{2n}{e}\Bigr)^{n/3}\exp\Bigl(\tfrac12(2n)^{2/3}-\tfrac16(2n)^{1/3}\Bigr).
\]
\end{itemize}
Comparing with the result of \cite{BGJ10a} recalled in the introduction, we obtain
$f^{(1)}_n/g^{(1)}_n\sim(C_1/C')\,n^{1/3}$. Figure~\ref{fig:p13} shows this growth, in contrast
with the convergence of $f^{(3)}_n/g^{(3)}_n$ established in Proposition~\ref{prop:olivierp}.
\erem

\begin{proof}[Proof of Theorem~\ref{thm:mainp}]
Part~(i) is Proposition~\ref{prop:growthp}, since $g^{(p)}_n\le f^{(p)}_n$. Part~(ii) is
Proposition~\ref{thm:fp}.
\end{proof}

\section{Concluding remarks}\label{sec:concl}

We have determined the asymptotic number of $BCK(p)$-terms, up to a constant factor, for every
$p\ge2$, and explicitly for $p=2$. The counting sequences exhibit fractional powers of the size $n$
in the exponent, and the Borel transform of the linearized generating function involves
exponentials of fractional powers of the variable. Other enumeration problems where powers of
$\sqrt n$ appear in the counting sequence are presented in \cite{Be74, BeCa78}. Fill et
al.~\cite{FJW12} presented an enumeration problem whose solution involves fractional powers which
are not multiples of $1/2$. Maybe that problem fits into the framework presented here; it might
then be possible to solve what the authors called a ``herculean task'', i.e. to obtain the full
asymptotic expansion of the considered sequence.

Several questions remain open.
\begin{itemize}
\item The constants $A$, $A_p$ and $\kappa_p=\lim f^{(p)}_n/g^{(p)}_n$ are not known in closed form.
\item Computer experiments for $p\le5$ suggest that $\sigma_p=-\frac{3p}{2(p+1)}$, that the
coefficient of $z$ in the exponent of Theorem~\ref{thm:Hp} equals $1$, and that the minimal
differential equation for $H_p$ has order $2p+1$ and coefficients of degree $p$. These facts are
not proved here.
\item For $p=1$, the quotient $f^{(1)}_n/g^{(1)}_n$ grows like $n^{1/3}$ (Remark~\ref{rem:p1}). Our
proof of this fact relies on the explicit result of \cite{BGJ10a}. A direct proof, along the lines
of Proposition~\ref{prop:olivierp}, would clarify the role of the quadratic term in this
borderline case.
\end{itemize}

\section{Concluding remarks}\label{sec:concl}

% ... (texte inchangé de la section 8) ...

\end{document}